\documentclass[xcolor=dvipsnames,svgnames,table,reqno,11pt]{amsart}

\usepackage[left=1.5in, right=1.5in, bottom=1.25in, height=8in]
{geometry}

\input xy
\xyoption{all}
\usepackage{accents}
\usepackage{epsfig}
\usepackage{xcolor}
\usepackage{amsthm}
\usepackage{amssymb}
\usepackage{bm} 
\usepackage{amsmath}
\usepackage{amscd}
\usepackage{amsopn}
\usepackage{graphicx}

\usepackage{tikz}
\usetikzlibrary{shapes,shadows,arrows}
\usepackage{tikz-cd}

\usepackage{xspace}

\usepackage{hhline}
\usepackage{easybmat}
\usepackage{caption}   
\usepackage{relsize}

\usepackage{url}
\usepackage{enumitem, hyperref}\hypersetup{colorlinks}

\usepackage{stmaryrd}

\usepackage[T1]{fontenc}

\newif\ifshowchanges
\showchangestrue

\newif\ifshowredchanges
\showredchangestrue

\colorlet{purpleB70}{blue!70!red}

\colorlet{orangeR65}{red!65!yellow}

\definecolor{red2}{HTML}{d41173}

\definecolor{neongreen}{HTML}{1bf702}

\definecolor{radicalred}{HTML}{FF355E}

\definecolor{denim}{HTML}{1560BD}

\definecolor{darkcyan}{rgb}{0.0, 0.55, 0.55}

\definecolor{cilek}{HTML}{FF43A4}

\definecolor{mor}{HTML}{9F00C5}

\definecolor{rorange}{RGB}{237,110,0} 

\definecolor{phlox}{rgb}{0.87, 0.0, 1.0}

\definecolor{fluorescentpink}{HTML}{FF1493}

\definecolor{napiergreen}{rgb}{0.16, 0.5, 0.0}

\definecolor{kellygreen}{rgb}{0.3, 0.73, 0.09}

\definecolor{parisgreen}{HTML}{ 50C878 }

\definecolor{palatinateblue}{rgb}{0.15, 0.23, 0.89}

\definecolor{ceruleanblue}{rgb}{0.16, 0.32, 0.75}

\definecolor{brandeisblue}{rgb}{0.0, 0.44, 1.0}

\definecolor{KLMblue}{HTML}{0FC0FC}

\definecolor{cinnamon}{rgb}{0.82, 0.41, 0.12}

\definecolor{darkorange}{rgb}{1.0, 0.55, 0.0}

\definecolor{darktangerine}{rgb}{1.0, 0.66, 0.07}

\definecolor{deepcarrotorange}{rgb}{0.91, 0.41, 0.17}

\definecolor{internationalorange}{HTML}{FF4F00}

\definecolor{persimmon}{HTML}{EC5800}

\definecolor{pumpkin}{HTML}{FF7518}

\definecolor{darkred}{rgb}{1,0,0} 
\definecolor{darkgreen}{rgb}{0,0.7,0}
\definecolor{darkblue}{rgb}{0,0,1}

\hypersetup{colorlinks,
linkcolor=darkblue,
filecolor=darkgreen,
urlcolor=darkred,
citecolor=darkgreen}

\makeatletter
\def\reflb#1#2{\begingroup
    #2%
    \def\@currentlabel{#2}%
    \phantomsection\label{#1}\endgroup
}
\makeatother

\numberwithin{equation}{section}
\newtheorem{Theorem}{Theorem}
\numberwithin{Theorem}{section}

\newtheorem{TheoremX}{Theorem}

\newtheorem{CorollaryX}[TheoremX]{Corollary}

\newtheorem   {Lemma}[Theorem]{Lemma}

\newtheorem*  {Conjecture}{Conjecture}
\newtheorem   {Proposition}[Theorem]{Proposition}
\newtheorem   {Corollary}[Theorem]{Corollary}
\theoremstyle {definition}
\newtheorem   {Definition}[Theorem]{Definition}
\theoremstyle {remark}
\newtheorem   {Remark}[Theorem]{Remark}

\def    \eps    {\epsilon}

\newcommand{\FF}{{\mathcal F}}

\newcommand{\CA}{{\mathcal A}}

\newcommand{\CI}{{\mathcal I}}

\newcommand{\CM}{{\mathcal M}}
\newcommand{\CN}{{\mathcal N}}
\newcommand{\CU}{{\mathcal U}}

\newcommand{\supp}{\operatorname{supp}}

\newcommand{\id}{{\mathit id}}

\newcommand{\const}{{\mathit const}}

\newcommand{\fe}{{\mathfrak e}}

\newcommand{\LL}{{\mathcal L}}

\newcommand{\PP}{{\mathcal P}}

\newcommand{\CW}{{\mathcal W}}

\def    \F      {{\mathbb F}}

\def    \C      {{\mathbb C}}
\def    \R      {{\mathbb R}}

\def    \Z      {{\mathbb Z}}
\def    \N      {{\mathbb N}}
\def    \Q      {{\mathbb Q}}

\def    \CP     {{\mathbb C}{\mathbb P}}

\def    \12     {{\frac{1}{2}}}

\def    \p      {\partial}

\def    \SH     {\operatorname{SH}}
\def    \COH    {\operatorname{CH}}
\def    \Sp     {\operatorname{Sp}}

\def    \tSp    {\operatorname{\widetilde{Sp}}}

\def    \HF     {\operatorname{HF}}

\def    \vDelta  {\vec{\Delta}}

\def    \mult  {\mathrm{h}}
\def    \lp {\mathit{lp}}

\def    \hmu   {\operatorname{\hat{\mu}}}

\def    \hh    {\mathrm{h}}

\newcommand \ee {\mathrm{e}}

\begin{document}


\setlength{\smallskipamount}{6pt}
\setlength{\medskipamount}{10pt}
\setlength{\bigskipamount}{16pt}





\title [On the Hofer--Zehnder conjecture for Reeb flows]{On the
  Hofer--Zehnder conjecture for Reeb flows without symmetry}

\author[Erman \c C\. inel\. i]{Erman \c C\. inel\. i}
\author[Viktor Ginzburg]{Viktor L. Ginzburg}
\author[Ba\c sak G\"urel]{Ba\c sak Z. G\"urel}

\address{E\c C: ETH Z\"urich, R\"amistrasse 101, 8092 Z\"urich,
  Switzerland}
\email{erman.cineli@math.ethz.ch}

\address{VG: Department of Mathematics, UC Santa Cruz, Santa Cruz, CA
  95064, USA} \email{ginzburg@ucsc.edu}

\address{BG: Kanouff School of Data, Mathematical and Statistical
  Sciences, University of Central Florida, Orlando, FL 32816, USA}
\email{basak.gurel@ucf.edu}



\date{\today} 

\thanks{The work is partially supported by the NSF grants DMS-2304207
  and DMS-2604231 (BG), DMS-2304206 and DMS-2604230 (VG), and the
  Simons Foundation grants SFI-MPS-TSM-00025241 (BG) and
  MP-TSM-00002529 (VG)}

\begin{abstract}
  This paper continues our investigation of the Hofer--Zehnder
  conjecture for dynamically convex, non-degenerate Reeb
  pseudo-rotations of $S^{2n-1}$. The main result of the paper is
  that, for $n\geq 3$, such a flow has at most $n-3$ alternating prime
  closed orbits and hence at most $2n-3$ prime closed orbits in
  total. In particular, the flow on $S^5$ has exactly three prime
  closed orbits, all of which are non-alternating and elliptic. In
  addition, we prove that a dynamically convex Reeb pseudo-rotation on
  a sphere has at least two non-degenerate, irrationally elliptic
  orbits.
  \end{abstract}

\maketitle

\vspace{-0.2in}


\tableofcontents

\section{Introduction and main results}
\label{sec:intro+results}

\subsection{Introduction}
\label{sec:intro}
In this paper, we continue our investigation, initiated in
\cite{CGG:HZ}, of the Hofer--Zehnder (HZ) conjecture for dynamically
convex, non-degenerate Reeb flows with finitely many prime closed
orbits on star-shaped hypersurfaces in $\R^{2n}$. We show that, for
$2n=6$, such a flow has exactly three prime closed orbits, settling
the conjecture for $S^5$ in the dynamically convex, non-degenerate
case.  In all dimensions $2n \geq6$, the number of orbits is bounded
from above by $2n-3$.  We also touch upon properties such as the
ellipticity of prime orbits in this setting. For instance, the orbits
are elliptic when $2n=6$.

For $2n\geq 6$, the central question in this area is the following
well-known conjecture, which goes back at least to \cite{Ek}.

\begin{Conjecture}[The $n$-or-$\infty$ conjecture]
  Let $M$ be the boundary of a star-shaped domain in $\R^{2n}$, which
  we assume to be smooth. Then

  \begin{itemize}

  \item[\reflb{C1}{\rm{(C-M)}}] {\rm Multiplicity or Ekeland's
      conjecture:} the Reeb flow on $M$ has at least $n$ prime closed
    orbits;
      
  \item[\reflb{C2}{\rm{(C-HZ)}}] {\rm HZ-conjecture:} the Reeb flow on
    $M$ has exactly $n$ prime closed orbits whenever the flow is a
    Reeb pseudo-rotation, i.e., the number of prime closed orbits is
    finite.

  \end{itemize}
\end{Conjecture}

Referring the reader to \cite{CGG:HZ} for a more detailed discussion,
here we only mention that very little is known about either part of
the conjecture in such generality. For instance, without extra
conditions on $M$, it is not even known that there must be more than
one prime orbit or more than two if the flow is
non-degenerate. However, once some additional restrictions on $M$ are
imposed, usually along the lines of convexity or dynamical convexity
or symmetry, much more is known.

One of the more natural assumptions of this type is dynamical
convexity: the requirement that all orbits have Conley--Zehnder index
$n+1$ or greater. (Convex hypersurfaces are dynamically convex, but
the latter condition is in general less restrictive than convexity;
see \cite{CE1, CE2, CGHi, DGRZ} and also \cite{AM, GM}.) For
non-degenerate flows, Ekeland's conjecture was proved in \cite{LZ} for
convex hypersurfaces and in \cite{GuK} in the dynamically convex
case. It is also known, essentially from \cite{LZ}, that 
non-degenerate dynamically convex pseudo-rotations have at least
$n$ non-alternating prime closed orbits. (A closed Reeb orbit $x$ is
non-alternating if the Conley--Zehnder indices of the iterates $x^k$
have the same parity.) Moreover, in this case the number of
non-alternating prime closed orbits is exactly $n$; see \cite{CGG:HZ}.

Furthermore, again in the non-degenerate setting, the dynamical
convexity condition has been significantly relaxed to a very minor
restriction on the index; see \cite{DLLW1, DLLW2, GGM} and references
therein.  Without the non-degeneracy requirement, Ekeland's conjecture
was proved in \cite{CGG:HZ} using the full form of dynamical
convexity; see \cite{LZ} and also \cite{DL, GG:LS, Wa1, Wa2, WHL} for
earlier results.

Much less is known about the HZ-conjecture, which is the main focus of
this paper and, to the best of our knowledge, was originally stated in
\cite{Gu:pr}. Here the main direct predecessor of this work is
\cite[Thm.\ B]{CGG:HZ} establishing the conjecture for non-degenerate
dynamically convex flows on centrally symmetric hypersurfaces; see
\cite{DQ, CGGM:Inv} for related results.

In our key result, Theorem \ref{thm:p<n}, we show that the number $q$
of alternating prime closed orbits, which could hypothetically exist
in addition to exactly $n$ non-alternating ones, is bounded above by
$n-3$ when $n \geq 3$. (Here we again assume dynamical convexity and
non-degeneracy and, of course, that the flow is a Reeb
pseudo-rotation.) As a consequence, the total number of orbits is
bounded from above by $2n-3$ when $n \geq 3$. When $n\leq 3$, all
prime closed orbits are non-alternating and there are exactly $n$ of
them; see Corollary~\ref{thm:6D} and also Remark \ref{rmk:4D}. This
result removes the symmetry assumption from \cite[Thm.\ B]{CGG:HZ} for
$n=3$. The approach to the question, which we develop in this paper,
is substantially new, and we outline it in the next section.

From a slightly different perspective, the $n$-or-$\infty$ conjecture
is a question about Reeb pseudo-rotations. In the context of
symplectic dynamics, the two competing definitions of a
pseudo-rotation are that this is a Hamiltonian system (a Hamiltonian
diffeomorphism or a Reeb flow) with finitely many periodic orbits or,
alternatively, with the minimal possible number of periodic orbits
(perhaps, counted with multiplicity); see \cite{GG:PR, Sh:HZ,
  CGGM:Inv}. The HZ-conjecture asserts that the two definitions are
equivalent, just as \cite{Sh:HZ} showed for Hamiltonian
diffeomorphisms of $\CP^n$ and some other manifolds; see also
\cite{AL, BX, CGG:HZ0}. Overall, perhaps expectedly so, there is a
considerable similarity between the properties of Hamiltonian and Reeb
pseudo-rotations; see, e.g., \cite{AK, LRS} vs.\ \cite{Ka} or
\cite{GG:PR} vs.\ \cite{CGGM:Inv} or \cite{Sh:HZ} vs.\ \cite{CGG:HZ}.
 
Reeb and Hamiltonian pseudo-rotations are of interest because they may
have quite complicated dynamics beyond periodic orbits; for instance,
they can be ergodic, \cite{AK, Ka, LRS}. However, in all known
examples the periodic orbit data of a pseudo-rotation (the linearized
return maps and the actions) is the same as that of a suitably chosen
$S^1$-action. In particular, every known Reeb pseudo-rotation on
$S^{2n-1}$ is dynamically convex and non-degenerate, and in fact has
the same periodic orbit data as an irrational ellipsoid. In Corollary
\ref{cor:ellipsoids}, we show that, under minor non-resonance
conditions, the periodic orbit data of a non-degenerate dynamically
convex Reeb pseudo-rotation in dimension 6 matches that of an
ellipsoid; cf.\ \cite {GG:RvsPR} for the Hamiltonian counterpart.

We also briefly touch upon the properties of the prime orbits. For
instance, we show that all three prime orbits in Corollary
\ref{thm:6D} are elliptic. Furthermore, we prove in Theorem
\ref{thm:elliptic} that every dynamically convex Reeb flow on
$S^{2n-1}$ with finitely many prime closed orbits has at least two
non-degenerate, irrationally elliptic orbits. We will further discuss
these results and some of the related previous work (e.g., \cite{AM,
  DQ, LLW, LO, LZ}) in Section \ref{sec:elliptic}.

\begin{Remark}
  The results of this paper are mainly of interest for $2n\geq 6$.
  The case of dimension $2n-1=3$ has been extensively studied and is
  well understood. According to the 2-or-infinity conjecture, a Reeb
  flow on a closed 3-manifold has either exactly two or infinitely
  many prime closed orbits. This conjecture is almost proved in
  general, and lens spaces are the only manifolds admitting Reeb
  pseudo-rotations; see \cite{CGH, CGHP, CGHHL, CDR}. Furthermore, the
  conjecture has been completely proved for $S^3$. (See also
  \cite{GHHM} for a different proof of \ref{C1} for $S^3$.)  Hence,
  both \ref{C1} and \ref{C2} parts of our $n$-or-$\infty$ conjecture
  hold in this case. However, one cannot expect this clear-cut
  dichotomy for all Reeb flows in higher dimensions.
\end{Remark}

\subsection{Main results and a discussion}
\label{sec:results}

On the dynamics side, the results of this paper break down into two
groups: those on the HZ-conjecture, which is our main focus, and a
series of complementary observations on the existence of elliptic
orbits of dynamically convex pseudo-rotations.

\subsubsection{HZ-conjecture}
Throughout the paper, all star-shaped and Liouville domains are
assumed to have smooth boundary. Our main result is the following.

\begin{TheoremX}
  \label{thm:p<n}
  Assume that the Reeb flow on the boundary of a star-shaped domain
  $W\subset \R^{2n}$ is non-degenerate, dynamically convex, and is a
  pseudo-rotation, i.e., it has finitely many prime closed
  orbits. Then the flow has exactly $n$ non-alternating prime closed
  orbits and $ q \leq \max\{0, n-3\}$ alternating ones.
\end{TheoremX}

Here, a prime closed orbit $z$ is said to be \emph{alternating} if the
parity of the Conley--Zehnder index $\mu\big(z^k\big)$ alternates as a
function of $k$ and \emph{non-alternating} if the parity is constant;
see Section \ref{sec:prelim}. The novel part of the theorem is the
upper bound $ q \leq \max \{0,n-3 \}$. The statement that there are
exactly $n$ non-alternating orbits is a part of \cite[Thm.\
A]{CGG:HZ}. This upper bound, in particular, implies that in the
setting of the theorem the total number of prime orbits is bounded by
$2n-3$ when $n\geq 3$.

\begin{CorollaryX}
  \label{thm:6D}
  Assume that the Reeb flow on the boundary of a star-shaped domain
  $W\subset \R^{6}$ is non-degenerate, dynamically convex and has
  finitely many prime closed orbits. Then the flow has exactly three
  prime orbits, and these orbits are elliptic.
\end{CorollaryX}

Indeed, $ q \leq 0$ when $n\leq 3$ by Theorem \ref{thm:p<n}, i.e., the
flow has no alternating orbits. When $\p W$ is centrally symmetric, a
similar result (sans ellipticity) holds in all dimensions; see
\cite[Thm.\ B]{CGG:HZ}.

A word is in order about the proof of Theorem \ref{thm:p<n}, which
follows a different path from other results of this type. Consider a
non-degenerate dynamically convex Reeb pseudo-rotation on
$S^{2n-1}$. Let $\Phi_1,\dotsc,\Phi_n$ be the linearized flows along
the non-alternating prime closed orbits and $\Psi_1,\dotsc,\Psi_{ q }$
the linearized flows along the alternating prime closed orbits. These
are elements of $\tSp(2m)$ with $m=n-1$. The results from
\cite{CGG:HZ} impose severe restrictions on these maps and their
iterates $\Phi_i^k$ and $\Psi_j^k$, which we recall in Theorem
\ref{thm:background}. For instance, all ratios of the mean indices are
irrational, the Conley--Zehnder indices of the iterates $\Phi_i^k$
with $k\in\N$ and $\Psi_j^k$ for odd $k$ partition the set
$\{n+1,\, n+3, \dotsc\}$, and such iterates with larger index have
larger mean index. We call a collection of maps $\{\Phi_i,\Psi_j\}$
meeting these conditions a formal pseudo-rotation; see Definition
\ref{def:FPR}. Now, the natural question is for what values of $n$ and
$ q $ formal pseudo-rotations exist. For instance, non-existence of
formal pseudo-rotations for all $ q >0$ would prove \ref{C2} in the
non-degenerate dynamically convex case. However, in this direct form
the question appears to be intractable.

To answer this question, we go one step further in formalizing the
setting. Namely, it turns out that the reciprocals of the mean
indices, $X_i=\hmu(\Phi_i)^{-1}$ and $Y_j=\hmu(\Psi_j)^{-1}$, must
meet certain combinatorial or number-theoretic requirements as a
consequence of the conditions on the maps. We call a collection
$\CU=\{X_i,Y_j\}$ of positive real numbers satisfying these
requirements an ultra-formal pseudo-rotation; see Section
\ref{sec:UFPR}. Thus every restriction on $n$ and $ q $ for
ultra-formal pseudo-rotations translates to a restriction for formal
pseudo-rotations and then to a restriction for
pseudo-rotations. Furthermore, we show that the resulting map from
formal to ultra-formal pseudo-rotations is onto; see Theorem
\ref{cor:mean}.

The main tool in the proof of this correspondence is a relation,
explored in Section \ref{sec:frequency_matching}, between generalized
Beatty-type sequences and their frequencies, which also has an
application to the properties of the Conley--Zehnder index.

The existence question for ultra-formal pseudo-rotations is already
much more tractable, and we show that $ q \leq n-3$ when $n\geq 3$;
see Theorem \ref{thm:numbers}. As a consequence, we obtain Theorem
\ref{thm:p<n}. However, rather unexpectedly for us, on the
combinatorial level, this result cannot be significantly pushed
further. Namely, for every $n>3$, there exists an ultra-formal
pseudo-rotation with $ q = 1$, as we show in Theorem
\ref{thm:numbers-ex}. Therefore, we also have formal pseudo-rotations
with $n>3$ and $ q = 1$. Moreover, there exist ultra-formal (and hence
formal) pseudo-rotations with arbitrarily large~$ q $.

\subsubsection{Elliptic orbits}
\label{sec:elliptic}
The provenance of the last statement of Corollary \ref{thm:6D} -- that
all three orbits are elliptic -- deserves a brief discussion. The
statement is an immediate consequence of Theorem~\ref{thm:p<n} and the
result from \cite{DQ} asserting that in dimension six the prime orbits
are elliptic under a minor index requirement relaxing dynamical
convexity, provided that there are exactly three prime orbits. The
statement also generalizes a theorem from \cite{LO} that the Reeb flow
on a strictly convex hypersurface in $\R^6$ has three elliptic closed
orbits, provided that the flow is non-degenerate and has only finitely
many prime orbits. For the sake of completeness, we include a short
proof of the ellipticity statement in Section \ref{sec:LA}, along with
the proof of Theorem \ref{thm:p<n} and the corollary, as an immediate
application of a result from \cite{CGGM:Inv}. The argument is
essentially identical to the one in \cite{DQ}, up to changes in the
notation and framework.

\begin{Remark}[Irrational ellipticity]
  An elliptic closed orbit is said to be \emph{irrationally elliptic}
  if the linearized return map is semisimple and the eigenvalues are
  not roots of unity. The latter condition is automatically satisfied
  for non-degenerate flows, but the former is not. At least two out of
  the three elliptic orbits in \cite{DQ,LO} are irrationally elliptic,
  and the same is true for the orbits in Corollary \ref{thm:6D}. (See
  the proof of Theorem \ref{thm:elliptic} for a direct proof of this
  fact.)
\end{Remark}

The assertion that the orbits are elliptic can be further refined
under a minor non-resonance condition. Namely, combining Corollary
\ref{thm:6D} with \cite[Cor.\ F]{CGG:HZ}, we see that the linearized
flows along the closed Reeb orbits precisely match those for an
irrational ellipsoid in $\R^6$.

To be more specific, denote by
$\vDelta:=\{\Delta_1<\Delta_2< \Delta_3\}\in\R^3_+$ the ordered set of
the mean indices of the three orbits from Corollary \ref{thm:6D}.  (By
\cite[Thm.\ E]{CGG:HZ}, these orbits have distinct mean indices; see
also Theorem \ref{thm:background}.) Consider the solid ellipsoid
$W_{\vDelta} \subset \C^3=\R^{6}$ given by
\[
  \pi \sum_j \frac{|z_j|^2}{\Delta_j}\leq 1
\]
and its boundary $M_{\vDelta}=\p W_{\vDelta}$. The prime closed orbits
on $M_{\vDelta}$ also have mean indices $\Delta_1<\Delta_2<
\Delta_3$. (This follows from the resonance relation
$\sum \Delta_j^{-1}=1/2$ from \cite{Vi:R}; see also, e.g.,
\cite{GiK}.)  Furthermore, up to scaling, $M=\p W$ and $M_{\vDelta}$
have the same action spectra. Indeed, let us scale $M$ so that
$\hmu(x)=\CA(x)$ for all closed orbits $x$ on $M$. This is possible
again by \cite[Thm.\ E]{CGG:HZ}.  Then for both $M$ and $M_{\vDelta}$
the action spectrum is simply the union of the sequences
$\Delta_j\N$. As an immediate consequence of \cite[Cor.\ F]{CGG:HZ},
we have

\begin{CorollaryX}[Ellipsoid comparison]
  \label{cor:ellipsoids}
  Assume that the Reeb flow on the boundary of a star-shaped domain
  $W\subset \R^{6}$ is non-degenerate, dynamically convex and has
  finitely many prime closed orbits.  In addition, let us require the
  following non-resonance condition: for every $j$ the mean index
  ratios $\pm\Delta_j/\Delta_i$, $i\neq j$, are distinct modulo $\Z$.
  Then the linearized flows along the corresponding closed orbits on
  $M$ and $M_{\vDelta}$ are equal up to conjugation as elements of
  $\tSp(4)$.
\end{CorollaryX}

Note that the non-resonance condition from the corollary need not be
satisfied even for irrational ellipsoids, although in this case the
statement is trivially correct: an ellipsoid matches itself. We do not
know if this condition is purely technical, which seems more likely,
or essential.

\begin{Remark}
  \label{rmk:4D}
  Analogues of Corollaries \ref{thm:6D} and \ref{cor:ellipsoids} hold
  for $2n=4$, with two prime orbits in place of three. However, in
  this case the result is true in a much more general setting and the
  two orbits are also known to be elliptic; see \cite{CGH, CGHP,
    CGHHL, CDR}. The non-resonance condition is void and the existence
  of a matching ellipsoid simply follows from the fact that
  $\Delta_1^{-1}+\Delta_2^{-1}=1/2$; see \cite{Vi:R}.
  \end{Remark}

  Returning to the existence of elliptic orbits, we have the
  following.

  \begin{TheoremX}
    \label{thm:elliptic}
    Assume that the Reeb flow on the boundary of a star-shaped domain
    $W\subset \R^{2n}$, $n\geq 2$, is dynamically convex and has
    finitely many prime closed orbits. Then the flow has at least two
    non-degenerate, irrationally elliptic prime closed orbits.
  \end{TheoremX}

  Theorem \ref{thm:elliptic}, proved in Section \ref{sec:dynamics}, is
  only tangentially related to the other results of the paper and the
  machinery developed here. Rather, it is essentially an immediate
  consequence, by a standard argument, of the proof of \cite[Thm.\
  A]{CGG:HZ}. The theorem generalizes and refines several previous
  results on the existence of elliptic closed orbits for Reeb
  pseudo-rotations. For instance, in the setting of the theorem for
  convex hypersurfaces, the existence of at least one elliptic orbit
  was proved in \cite{LZ}. The result was extended to the dynamically
  convex case in \cite{AM}. The existence of at least two
  non-degenerate, irrationally elliptic closed orbits was established
  in \cite{LLW} under the additional assumption that there are exactly
  $n$ prime closed orbits. Finally, the existence of at least two
  irrationally elliptic closed orbits for non-degenerate Reeb
  pseudo-rotations on convex hypersurfaces was proved in \cite{LO}.
  We refer the reader to, e.g., \cite{DQ} for a detailed discussion of
  other results on elliptic orbits of Reeb pseudo-rotations and
  further references. (We do not touch here on the relevant but
  different question of the existence of non-hyperbolic orbits; see,
  e.g., \cite{AM22, AM26, DQ, GGM} and references therein.)
  
  \begin{Remark}
    \label{rmk:elliptic}
    The non-degeneracy of the orbits found in Theorem
    \ref{thm:elliptic} is essential to guarantee irrational
    ellipticity. Moreover, without non-degeneracy, the very notion of
    ellipticity (only unit eigenvalues) should be taken with a grain
    of salt. For example, on a hypersurface in $\R^4$, consider an
    orbit such that its return map is a standard monkey saddle. This
    orbit is totally degenerate and hence elliptic. However, it is
    locally maximal as an invariant set and, from this perspective,
    more similar to a hyperbolic orbit. For instance, the crossing
    energy bounds from \cite{ CGGM:Entropy, CGGM:Inv} apply to this
    orbit. Finally, let us point out that this class of questions
    stems from the general conjecture that the Reeb flow on the
    boundary of a star-shaped domain in $\R^{2n}$ has an elliptic
    closed orbit. Without going into a discussion of this conjecture
    and its history, we only mention here that the conjecture remains
    wide open, even for non-degenerate flows on convex hypersurfaces
    in $\R^4$. (See, e.g., \cite{AM:elliptic} for a more detailed
    account.) The condition that the flow is a pseudo-rotation
    completely changes the nature of the question and makes it much
    more tractable.
  \end{Remark}  

  \smallskip\noindent\subsection*{Acknowledgments and use of AI tools}
  The authors are grateful to Leonardo Macarini for his patient and
  detailed answers to their questions.

  The main idea and approach to the question are entirely due to the
  authors and build on their previous work. In carrying out this work,
  the authors used Claude Opus 5.5 and, to a lesser extent, Gemini 3.1
  Pro to assist in generating proofs through an iterative process
  directed by the authors who formulated the problems, posed the
  intermediate questions, and decided which directions to pursue. The
  authors wrote the paper, checked all the arguments in detail, and
  take full responsibility for its content.

\section{Preliminaries}
\label{sec:prelim}

\subsection{Conventions, notation and basics: symplectic linear
  algebra}
\label{sec:conventions}
In this section, we briefly review the relevant notions from
symplectic linear algebra and set notation and conventions.

Let $\Phi\in\tSp(2m)$ be a path parametrized by $[0,1]$, starting at
the identity and taken up to homotopy with fixed endpoints. \emph{We
  will refer to the eigenvalues of $\Phi(1)$ as the eigenvalues of
  $\Phi$ and call $\Phi$ a map}. Recall that $\Phi$ is said to be
\emph{non-degenerate} if one is not an eigenvalue and \emph{strongly
  non-degenerate} if all iterates $\Phi^k$, $k\in\N$, are
non-degenerate. A map $\Phi$ is \emph{elliptic} when all eigenvalues
are unit. Throughout the paper, we denote by $\hmu(\Phi)$ the
\emph{mean index} of $\Phi$ and, when $\Phi$ is non-degenerate, by
$\mu(\Phi)$ the \emph{Conley--Zehnder index}, referring the reader to,
e.g., \cite{SZ} for the definitions. Alternatively, the reader can
take \eqref{eq:mult-index} and \eqref{eq:mean} as the definitions of
these indices. Then, as for instance readily follows from these
formulas,
\begin{equation}
  \label{eq:mu-hmu}
|\mu(\Phi)-\hmu(\Phi)|<m .
\end{equation}
The mean index is homogeneous:
$\hmu\big(\Phi^k\big)=k\hmu(\Phi)$. (The Conley--Zehnder index is
not.) Then, by \eqref{eq:mu-hmu},
$\mu\big(\Phi^k\big)/k\to \hmu(\Phi)$ when $\Phi$ is strongly
non-degenerate.

The entries in the sequence of iterated indices $\mu\big(\Phi^k\big)$,
$k\in\N$, either all have the same parity or the parity alternates,
provided of course that $\Phi$ is strongly non-degenerate. We say that
$\Phi$ is \emph{non-alternating} in the former case and
\emph{alternating} in the latter. A map $\Phi$ is non-alternating if
and only if the number of real eigenvalues in the range $(-1, 0)$ is
even and hence alternating if and only if this number is odd.

Our next goal is to associate to $\Phi$ several numerical invariants
which determine its index behavior under iterations. For the sake of
simplicity, let us first assume that $\Phi$ is semisimple, i.e., for
every eigenvalue its algebraic multiplicity is equal to the geometric
multiplicity, although we do not require the eigenvalues to be
distinct. Then it is easy to see that we can write $\Phi$ as a product
of a loop $\psi$ with the direct sum
$\Phi_{h}\oplus \Phi_{-h}\oplus \Phi_e$, where the maps $\Phi_{\pm h}$
and $\Phi_e$ are as follows.

The map $\Phi_h$ is hyperbolic with complex or positive real
eigenvalues and zero mean index. An equivalent condition is that
$\Phi_h(t)$ is hyperbolic with complex or positive real eigenvalues
for all $t\in (0,\,1]$.  The second term $\Phi_{-h}$ is hyperbolic
with negative real eigenvalues. We denote by $\hh_{-}(\Phi)$ the
number of negative hyperbolic blocks, i.e., half of the dimension of
the domain of $\Phi_{-h}$.  As an element of the universal covering,
$\Phi_{-h}$ is specified by the condition that
$\hmu(\Phi_{-h})=\hh_{-}(\Phi)$ or, equivalently, by the condition
that it is connected to $I$ by the counterclockwise rotation by $\pi$
and a hyperbolic transformation.  Thus $\Phi$ is alternating exactly
when $\hh_{-}(\Phi)$ is odd. Finally, let $\hh(\Phi)$ be half of the
dimension of the domain of $\Phi_{h}\oplus \Phi_{-h}$.

The remaining term $\Phi_e$ is elliptic and, by our assumption,
semisimple. Up to conjugation, it decomposes as a direct sum of
counterclockwise rotations $R_{\theta_i}$ by angles
$2\pi\theta_i\in (0,2\pi)$. We denote by $\ee(\Phi)$ half of the
dimension of the domain of $\Phi_e$. Note that
$\hh(\Phi)+\ee(\Phi)=m$. We call the ``logarithmic eigenvalues''
$\theta_i$ the \emph{rotation numbers} of $\Phi$ and write
$\theta_i(\Phi)$ when it is necessary to emphasize the dependence on
$\Phi$.  (The ordering of this multiset is not essential for our
purposes, but one can use the non-decreasing order to make literal
sense of this expression.) We will denote the multiset of the rotation
numbers by $\Theta(\Phi)$. This normalization of $\Phi_e$ is not
standard and is different from, for instance, the one used in
\cite{CGG:HZ, GG:LS, GG:RvsPR}. However, it is convenient for our
purposes as it eliminates negative terms in \eqref{eq:mult-index} and
\eqref{eq:mean} other than $2\lp(\Phi)$.

The loop part $\psi$ is then determined by $\Phi$ and the above
normalizations of $\Phi_{\pm h}$ and $\Phi_e$ viewed as elements of
the universal covering of the group of linear symplectic
transformations.

Set
$$
\lp(\Phi):=\hmu(\psi)/2.
$$
In other words, $2\lp(\Phi)$ is the Maslov index of $\psi$, i.e.,
$\lp(\Phi)$ is the image of the free homotopy class of $\psi$ under
the isomorphism $\pi_1\big(\Sp(2m)\big)\cong\Z$ sending the
counterclockwise rotation by $2\pi$ in a symplectic plane to one.

It is easy to extend these definitions to non-degenerate maps $\Phi$
which are not semisimple. To this end, we may, for instance, connect
$\Phi$ to a semisimple map $\Phi'$ with the same spectrum as $\Phi$
via a family of isospectral transformations, i.e., a family of
transformations with constant spectrum; see, e.g., \cite[Lemma
5.3]{CGG:HZ}. We set $\mult_{-}(\Phi):=\mult_{-}(\Phi')$,
$\lp(\Phi):=\lp(\Phi')$, $\ee(\Phi):=\ee(\Phi')$ and
$\Theta(\Phi):=\Theta(\Phi')$. Then it is routine to show that the
resulting invariants are independent of the choice of $\Phi'$.  It is
also clear that the invariants are additive with respect to direct
sum.

Essentially by definition,
\begin{equation}
\label{eq:mult-index}
\mu(\Phi)=2\lp(\Phi)+\mult_{-}(\Phi)+\ee(\Phi)
\end{equation}
and
\begin{equation}
\label{eq:mean}
\hmu(\Phi)=2\lp(\Phi)+\hh_{-}(\Phi)+
2\sum_{i=1}^{\ee(\Phi)}\theta_i(\Phi).
\end{equation}
Therefore, 
$$
|\mu(\Phi)-\hmu(\Phi)| \leq \ee(\Phi),
$$
with strict inequality when $\ee(\Phi) >0$, a refinement of
\eqref{eq:mu-hmu}. Finally, when $\Phi$ is strongly non-degenerate,
for all $k\in\N$, we have
\begin{equation}
\label{eq:mu-k}
\mu\big(\Phi^k\big)=k\big(2\lp(\Phi)+
\hh_{-}(\Phi)\big)+\ee(\Phi)
+2\sum_{i=1}^{\ee(\Phi)}\lfloor k\theta_i(\Phi)\rfloor .
\end{equation}

\subsection{Background results from symplectic dynamics}
\label{sec:dynamics}
Consider now the Reeb flow on the boundary $S^{2n-1}=\p W$ of a
star-shaped domain $W\subset \R^{2n}$. For every periodic orbit $x$,
pick a trivialization of the contact structure $\xi$ on $S^{2n-1}$
along $x$ so that the induced trivialization of the top complex wedge
$\wedge_{\C}^{n-1}\xi|_x$ extends to $S^{2n-1}$. Let
$\Phi\in \tSp\big(2(n-1)\big)$ be the linearized Reeb flow along
$x$. This map, as an element of the universal covering, is independent
of the choice of the trivialization, and the linear algebra concepts
about $\Phi$ translate to $x$. For instance, $x$ is non-degenerate if
$\Phi$ is non-degenerate; the Conley--Zehnder index $\mu(x)$ of $x$ is
$\mu(\Phi)$; $x$ is (non-)alternating if so is $\Phi$; etc.  A
non-degenerate closed orbit is called \emph{bad} if it is an even
iterate of an alternating orbit; all other closed orbits are
\emph{good}.

Furthermore, the flow is \emph{dynamically convex} if $\mu(x)\geq n+1$
for all closed Reeb orbits $x$. This condition, introduced in
\cite{HWZ}, follows from convexity, but there exist star-shaped
domains which are dynamically convex but not symplectomorphic to
convex domains; see \cite{CE1, CE2} and also \cite{CGHi, DGRZ}. We
refer the reader to, e.g., \cite{GG:LS} for a discussion of dynamical
convexity in the context of symplectic linear algebra.

The following theorem is the key input from symplectic dynamics on
which our proofs rely. It comprises parts of \cite[Thms.\ A and E,
Cor.\ D]{CGG:HZ}.

\begin{Theorem}[\cite{CGG:HZ}]
  \label{thm:background}
  Assume that the Reeb flow on the boundary of a star-shaped domain
  $W\subset \R^{2n}$ is non-degenerate, dynamically convex and has
  finitely many prime closed orbits. Then
  \begin{itemize}
   
  \item the flow has exactly $n$ non-alternating prime closed orbits;

  \item the Conley--Zehnder indices of good iterates partition the set
    $\{n+1, n+3, \ldots\}$;

  \item all closed orbits have the same mean index-action ratio;

  \item all pairwise ratios of the mean indices of the distinct prime
    closed orbits are irrational, and, in particular, all closed
    orbits have distinct mean indices;
    
  \item for any two closed orbits $z$ and $z'$ (good or bad),
    $\mu(z)\geq\mu(z')$ whenever $\hmu(z)>\hmu(z')$, and the former
    inequality is strict when both orbits are good.
    
  \end{itemize}
\end{Theorem}

\begin{Remark}
\label{rmk:set-vs-multiset}
Note that, by the penultimate point of the theorem, the linearized
flows along prime closed orbits are distinct. Hence, these flows form
a genuine subset of $\tSp(2m)$, $m=n-1$, rather than a multiset.
\end{Remark}

Referring the reader to \cite{CGG:HZ} for detailed proofs, here we
only point out the key steps in the argument. Denote by $\SH^t(W)$ the
filtered symplectic homology of $W$ with coefficients in a fixed field
$\F$, and set $D_t=\dim \SH^t(W)$. The central point of the proof is
that $D_t=1$ for all $t>0$. (This identity does not require
non-degeneracy.) Then, for instance, the third item is a consequence
of this fact and a result from \cite{GG:LS}. The first ingredient of
the proof is an \emph{a priori} bound on the bar length in the
persistence module $\SH^t(W)$: all bars have finite length, which is
bounded above by some constant independent of $t$ and $\F$; see
\cite{GS}. Then, $D_{\max}:=\max_{\tau} D_\tau<\infty$ due to this
upper bound and the fact that $\hmu(x)>0$ for all closed Reeb orbits
$x$ by the dynamical convexity condition.  The fact that $D_t\geq 1$
for all $\F$, including $\F=\Q$, follows from the Smith inequality,
$D_{pt}\geq D_t$, for $\F=\F_p$ for all primes $p>0$; see
\cite{Se,ShZ}. (These two ingredients appear to be common to all
results of this type starting with \cite{Sh:HZ}.)  Furthermore, by the
Smith inequality again, $D_t=D_{\max}$ for $t$ in a sequence of
arbitrarily long action intervals $I$. Finally, we show that
$D_t\leq 1$ for some $I$ and $t\in I$, and hence
$1\leq D_t\leq D_{\max}\leq 1$ for $t>0$, i.e., $D_t=1$. This is done
via the index recurrence theorem, which is essentially an input from
the geometry of numbers, heavily relying on the dynamical convexity;
see \cite[Sec.\ 3.2]{CGG:HZ}.

We conclude this section by sketching how Theorem \ref{thm:elliptic}
follows from the proof of \cite[Thm.\ A]{CGG:HZ}.

\begin{proof}[Proof of Theorem \ref{thm:elliptic}]
  The proof of \cite[Thm.\ A]{CGG:HZ} entails, for every $\eps>0$, the
  construction of an interval
  $$
  \CI=[d-n+1,d+n-1],
  $$
  where $d$ is an even integer, $n$ distinct prime closed orbits
  $x_1,\ldots, x_n$, and positive even integers $k_1,\ldots,k_n$
  satisfying several conditions. (See \cite[Sec.\ 4.1.1]{CGG:HZ}.)
  These conditions are
\begin{itemize}
\item $\big|\hmu\big(x_i^{k_i}\big)-d\big|<\eps$ for all $i$;
\item the local $S^1$-equivariant symplectic homology
  $\COH\big(x_i^{k_i}\big)$ of $x_i^{k_i}$ over $\Q$ is
  one-dimensional and hence concentrated in one degree, which we
  denote by $\deg\big(x_i^{k_i}\big)$;
\item
  $\big|\hmu\big(x_i^{k_i}\big)-\deg\big(x_i^{k_i}\big) \big|\leq
  n-1$;
\item for every integer $\ell\in \CI$ of the same parity as $n-1$
  there exists exactly one $x_i$ with $\deg\big(x_i^{k_i}\big)=\ell$.
\end{itemize}

Next, recall that whenever a closed (possibly degenerate) orbit $y$ is
not elliptic, we have
\begin{equation}
  \label{eq:non-elliptic}
  \supp \COH(y)\subset [\hmu(y)-n+2,\hmu(y)+n-2].
\end{equation}  
Here $\supp \COH(y)$ is the range of degrees where the local homology
is non-zero; see, e.g., \cite[Sec.\ 2.4.1]{CGG:HZ} or \cite[Sec.\
2.3.1]{GG:LS}.

It follows that $\deg\big(x_i^{k_i}\big)$ is an interior point of
$\CI$ whenever $x_i^{k_i}$ is not elliptic. Hence, in this case, and
when $\eps>0$ is small enough, we have
\begin{equation}
\label{eq:non-elliptic2}
\big|\deg\big(x_i^{k_i}\big)-d\big|<n-1.
\end{equation}

Therefore, we must have two elliptic orbits, say $x_{i_\pm}$, with
degrees at the endpoints of $\CI$. Finally, a prime orbit is elliptic
if and only if any one of its iterates is elliptic (and if and only if
all its iterates are elliptic). Hence, both orbits $x_{i_\pm}$ are
elliptic, although we have not yet shown that these orbits are
non-degenerate.

To establish the non-degeneracy of these orbits, set
$y_{\pm}=x_{i_{\pm}}^{k_{i_{\pm}}}$. The first step is to observe that
the local Floer homology $\HF(\varphi_{\pm})$ of the return map
$\varphi_{\pm}$ of $y_{\pm}$ is supported in one degree, and hence
\cite[Lemma 4.4]{CGG:HZ} applies to the admissible odd iterates of
$y_{\pm}$. Namely, since $\deg(y_{\pm})= d \pm (n-1)$, we have
$\HF_{d \pm (n-1)}(\varphi_{\pm}) \neq 0$. In other words, the local
homology is supported in the top/bottom degree, respectively. Then, by
\cite{GG:gap} and \cite[Lemma 3.6]{Gi:CC}, generating functions for
the degenerate parts of $\varphi_{\pm}$ have an isolated
maximum/minimum. As a consequence, we infer, in particular, that
$\dim \HF(\varphi_{\pm}) =1$. Next, as in the proof of \cite[Thm.\
A]{CGG:HZ}, we pick an ``index recurrence event'' with center $d'$ and
such that the first $N$ visible orbits (i.e., closed orbits with
non-zero local homology) on both sides of the event are odd admissible
iterates of the first $N$ visible orbits. We choose $N$ large enough
so that $y_{\pm}$ are among these orbits. Now, we have two copies of
$y_{\pm}$ on the two sides of the event and their degrees are
$d' \pm \deg(y_{\pm})$; see \cite[Thm.\ E]{CGG:HZ}. On the other hand,
by index recurrence and \cite[Lemma 4.4]{CGG:HZ}, the difference of
the two, $2\deg(y_{\pm})$, is equal to two times the Conley--Zehnder
index of the non-degenerate part of the linearized flow along
$y_{\pm}$. (Here we are assuming that the error in index recurrence,
i.e., the difference between $d'$ and the mean indices, is
sufficiently small.) For the orbits $y_{\pm}$ with degree at the
endpoints of $\CI$, this is possible only when the orbits are
non-degenerate; for otherwise we would have
$| d-\deg(y_{\pm}) | <n-1$.

It remains to show that the orbits $x_{i_\pm}$ are irrationally
elliptic. The assertion that no eigenvalue is a root of unity follows
from the non-degeneracy of the iterates
$y_{\pm}=x_{i_\pm}^{k_{i_\pm}}$ together with the fact that, by
construction, $k_{i_\pm}$ is divisible by the orders of the roots of
unity among the eigenvalues; see the proof of \cite[Thm.\
A]{CGG:HZ}. Next, it is not hard to infer from Williamson's normal
forms, \cite{Wi1, Wi2, Wi3}, that \eqref{eq:non-elliptic} also holds
when the return map of $y$ is non-degenerate but not semisimple. In
fact, then $|\mu(y)-\hmu(y)|\leq n-3$, provided that $-1$ is not an
eigenvalue. (Ultimately, this is a consequence of the fact that a
positive (negative) definite quadratic Hamiltonian is automatically
semisimple. Thus, non-degenerate non-semisimple blocks force a
cancellation of the first Krein-type eigenvalues.)  Applying this
observation to $y_{\pm}$, we see that \eqref{eq:non-elliptic2} would
also hold for these orbits if they were not semisimple, contradicting
$\deg(y_{\pm})= d \pm (n-1)$.  Therefore, the orbits $y_{\pm}$ and
hence $x_{i_\pm}$ are semisimple, and thus irrationally elliptic. This
completes the proof of the theorem.
\end{proof}

\section{Formal pseudo-rotations}
\label{sec:LA}
Our goal in this section is to formalize the statement of Theorem
\ref{thm:background} in the context of symplectic linear algebra. To
this end, we first need to introduce some terminology. Let
$\FF\subset \tSp(2m)$ be a countable collection of non-degenerate
symplectic linear maps or, to be more precise, 
elements of the universal covering of the linear symplectic group. We
say that $\FF$ \emph{LS-partitions} a set $K\subset \Z$ if all
elements of $\FF$ have distinct Conley--Zehnder indices, $K$ is
exactly the set of such indices, and the elements with larger
Conley--Zehnder index have a larger mean index, i.e., $\mu(F)>\mu(F')$
iff $\hmu(F)>\hmu(F')$ for any two elements $F$ and $F'$ of
$\FF$. (Here LS stands for Lusternik--Schnirelmann.) In particular,
the elements of $\FF$ have distinct mean indices.  Furthermore, an
iterate $F^k$, $k\in\N$, of a strongly non-degenerate map
$F\in \tSp(2m)$ is said to be 
\emph{good} if $k$ is odd and $F$ is alternating or $k$ is arbitrary
and $F$ is non-alternating.

Finally, it will be convenient in the rest of the paper to adopt the
following definition.

\begin{Definition}
  \label{def:FPR}
  Let $\Phi_1, \ldots, \Phi_{m+1}$ be distinct strongly non-degenerate
  non-alternating elements of $\tSp(2m)$ and let
  $\Psi_1, \ldots, \Psi_{ q }$ be a finite set (possibly empty) of
  strongly non-degenerate alternating elements. We say that the
  collection of elements $\{\Phi_i,\Psi_j\}$ is a \emph{formal
    pseudo-rotation} if the following three conditions are met:
  \begin{itemize}

  \item[\reflb{cond-irr}{\rm{(IR)}}] all pairwise ratios of the mean
    indices $\hmu(\Phi_i)$ and $\hmu(\Psi_j)$ are irrational, except
    of course the ratios
    $\hmu(\Phi_i)/\hmu(\Phi_i)=1=\hmu(\Psi_j)/\hmu(\Psi_j)$;

  \item[\reflb{cond-monotone}{\rm{(MO)}}]
    $\mu(\varphi)\geq\mu(\varphi')$ for any two iterates $\varphi$ and
    $\varphi'$ (good or bad) with $\hmu(\varphi)>\hmu(\varphi')$;

  \item[\reflb{cond-LSpart}{\rm{(LS)}}] the good iterates $\varphi$ of
    $\Phi_i$ and $\Psi_j$ LS-partition the set
    $K:=\{m+2, m+4, \ldots\}\subset \Z$.
    
\end{itemize}
\end{Definition}

In what follows, it is essential to separate the sets $\{\Phi_i\}$ and
$\{\Psi_j\}$ and treat them differently. Thus, in the definition, a
formal pseudo-rotation is actually a pair of sets.

Our next theorem, established in Sections \ref{sec:reduction} and
\ref{sec:proofs}, is the combinatorial underpinning of Theorem
\ref{thm:p<n}.

\begin{Theorem}
  \label{thm:LA-p_leq_m}
  Let $\Phi_1, \ldots, \Phi_{m+1}$ and $\Psi_1, \ldots, \Psi_{ q }$ be
  a formal pseudo-rotation in $\tSp(2m)$. Then
  $ q \leq \max\{0,m-2\}$.
\end{Theorem}

\begin{Remark}
  Condition \ref{cond-irr} can be rephrased as the condition that all
  iterates (good or bad) have distinct mean indices. Condition
  \ref{cond-monotone} is equivalent to the condition that for all
  iterates the index is a monotone (but possibly not strictly
  monotone) function of the mean index. Finally, the LS-partition
  requirement, \ref{cond-LSpart}, asserts that $K$ is in one-to-one
  correspondence with the set of good iterates and the mean index is a
  strictly monotone increasing function on $K$, or, equivalently, that
  the index is a strictly monotone increasing function of the mean
  index for good iterates.

  Note that \ref{cond-irr} and \ref{cond-monotone} almost, but not
  quite, follow from \ref{cond-LSpart}. In fact, they would if all
  iterates were good.  However, with the restriction to only odd
  iterates for alternating orbits, we can have an LS-partition with a
  rational ratio involving an alternating orbit. Moreover, two bad
  iterates with consecutive mean indices have the same Conley--Zehnder
  index, i.e., in this case the index is a monotone, but not strictly
  monotone, function of the mean index.
  \end{Remark}

\begin{Corollary}
  \label{thm:LA-m=2}
  In the setting of Theorem \ref{thm:LA-p_leq_m}, assume that
  $m\leq 2$. Then $ q =0$, i.e., the collection of alternating
  elements is empty.
\end{Corollary}

\begin{proof}[Proof of Theorem \ref{thm:p<n} and Corollary \ref{thm:6D}]
  In the setting of the theorems, consider the linearized flow in
  $S^{2n-1}=\p W$ along each prime closed orbit $x$. (Recall that this
  flow depends on the trivialization of the contact structure on
  $S^{2n-1}$ along $x$. The trivialization is chosen so that the
  induced trivialization of the top complex wedge
  $\wedge_{\C}^{n-1}\xi|_x$ extends to $S^{2n-1}$.) By definition,
  non-alternating periodic orbits give rise to non-alternating maps
  $\Phi_i\in \tSp(2m)$, where $m=n-1$, and the alternating orbits give
  rise to alternating maps $\Psi_j$. By Theorem \ref{thm:background},
  the conditions of Theorem \ref{thm:LA-p_leq_m} and Corollary
  \ref{thm:LA-m=2} are met, and Theorem \ref{thm:p<n} and Corollary
  \ref{thm:6D} follow with the exception of the statement that the
  three orbits in the corollary are elliptic; cf.\ \cite{DQ}.

  To see this, note that these orbits are non-alternating and have
  even index. A non-degenerate, non-alternating element of $\tSp(4)$
  can have an even index only when it is either elliptic or
  hyperbolic. None of the prime orbits can be hyperbolic since a
  dynamically convex Reeb flow on $S^{2n-1}$ with a prime hyperbolic
  orbit has infinitely many prime orbits; see \cite{CGGM:Inv}. Hence,
  all three closed orbits are elliptic.
 \end{proof}

 Our next result, proved in Section \ref{sec:pf-numbers-ex}, shows
 that Corollary \ref{thm:LA-m=2}, i.e., the statement that $ q =0$
 when $m\leq 2$, does not extend to dimensions $2m\geq 6$. This fact,
 which we find quite surprising, sets a natural limit to the
 combinatorial approach to the proof of the HZ-conjecture.

\begin{Theorem}
  \label{thm:LA-ex}
  For every integer $m\geq 3$, there exists a formal pseudo-rotation
  with $ q = 1$. Moreover, there exist formal pseudo-rotations with
  arbitrarily large $ q $ (and then, of course, large $m$).
\end{Theorem}

\section{Frequency matching}
\label{sec:frequency_matching}

In this section, we state a technical result on generalized
Beatty-type sequences. This is Proposition \ref{lemma:FM}, proved in
Section \ref{sec:pf-lemma-FM}. In essence, the statement is that a
Beatty-type sequence determines the frequencies up to natural
cancellations. The result plays a key role in the proof of Theorem
\ref{thm:LA-p_leq_m} and, in particular, in the proof of Theorem
\ref{cor:mean}. Furthermore, it has applications to the properties of
the Conley--Zehnder index (Theorem \ref{thm:index-spec}), and might
also be of independent interest.

\subsection{Detecting frequencies}
\label{sec:detect_frequencies}

Consider two finite multisets $A=\{\alpha_1,\ldots,\alpha_{|A|}\}$ and
$B=\{\beta_1,\ldots,\beta_{|B|}\}$ of positive irrational numbers. We
say that $\alpha_i$ and $\beta_j$ are \emph{externally matched} if
$\alpha_i=\beta_j\mod 1$. A pair $\alpha_i$ and $\alpha_j$ (or
$\beta_i$ and $\beta_j$) is \emph{internally matched} if
$\alpha_i+\alpha_j=0\mod 1$ (respectively, $\beta_i+\beta_j=0\mod
1$). Set
  \begin{equation}
    \label{eq:f_A}
  f_A(k)=\sum_i \{k\alpha_i\},
\end{equation}
where $k\in \N$ and $\{\cdot\}$ stands for the fractional part, and
likewise
  \begin{equation}
    \label{eq:f_B}
  f_B(k)=\sum_j \{k\beta_j\}.
\end{equation}

\begin{Proposition}[Frequency Matching]
\label{lemma:FM}
The sequence $f_A-f_B$ is constant if and only if the multiset union
$A\sqcup B$ breaks down into a union (perhaps, non-uniquely) of
internally or externally matched pairs. Furthermore, in this case,
$f_A-f_B=(|A|-|B|)/2$, and this is an integer.
\end{Proposition}

As an immediate consequence of the proposition we have the following.

\begin{Corollary}
 \label{cor:FM0}
 The sequences $f_A$ and $f_B$ are equal, i.e., $f_A(k)=f_B(k)$ for
 all $k\in \N$, if and only if $|A|=|B|$ and the multiset union
 $A\sqcup B$ breaks down into a union (perhaps, non-uniquely) of
 internally or externally matched pairs.
\end{Corollary}

Proposition \ref{lemma:FM} can be rephrased in terms of generalized
Beatty sequences. Namely, let $A$ and $B$ be as above. Set
  $$
  b_A(k)=\sum_i \lfloor k\alpha_i\rfloor,
  $$
  where $k\in \N$, and 
  $$
  b_B(k)=\sum_j \lfloor k\beta_j\rfloor.
  $$

  \begin{Corollary}[Frequency Matching in Beatty sequences]
\label{cor:FM}
We have $b_A(k)-b_B(k)=a k-c$ for all $k\in\N$, with $a$ and $c$
constant, if and only if the multiset $A\sqcup B$ breaks down into a
union of internally or externally matched pairs. In this case,
$a=\sum_{\alpha\in A}\alpha-\sum_{\beta\in B}\beta$ and
$c=(|A|-|B|)/2$ are both integers.
\end{Corollary} 

\begin{proof} The corollary follows immediately from Proposition
  \ref{lemma:FM} and the fact that
  $\lfloor k\alpha \rfloor =k\alpha-\{k\alpha\}$ for all positive
  irrational $\alpha$. To see that $a$ is an integer, it suffices to
  set $k=1$. Then $a=b_A(1)-b_B(1)+c$ with all integer terms on the
  right.
\end{proof}

It is not hard to see that Corollary \ref{cor:FM} is in fact
equivalent to Proposition \ref{lemma:FM}. As a particular case of the
corollary, we have

  \begin{Corollary}  
  \label{cor:FM2}
  The sequences $b_A$ and $b_B$ are equal if and only if $|A|=|B|$,
  the multiset $A\sqcup B$ breaks down into a union of internally or
  externally matched pairs, and $\sum_i\alpha_i=\sum_j\beta_j$.
\end{Corollary}

These results can be established as a consequence of \cite[Thm.\
1.8]{He}. For the sake of completeness, we include a direct and
relatively short proof relying on the machinery of almost periodic
functions in the sense of Besicovitch.

\subsection{Spectrum and the iterated Conley--Zehnder index}
\label{sec:Appl_CZ}
While Proposition \ref{lemma:FM} is central to the proofs of the main
results of the paper, it has other applications. In this section, we
touch upon a combinatorial consequence of the proposition and
Corollary \ref{cor:FM}, asserting, roughly speaking, that the sequence
of iterated Conley--Zehnder indices of a map determines its elliptic
eigenvalues and the size of its hyperbolic part, up to ambiguity
coming from cancellation of ``opposite eigenvalues''. This result was
originally established as a particular case of a more general theorem,
\cite[Thm.\ 4.2]{GG:RvsPR}, proved by a completely different ``brute
force'' combinatorial method.  (The reader should be aware that the
conventions used in \cite{GG:RvsPR} are different from those adopted
here.)  Proposition \ref{lemma:FM} and its corollaries provide a
``civilized'' proof of Theorem \ref{thm:index-spec}.

\begin{Theorem} [\cite{GG:RvsPR}]
\label{thm:index-spec}
Let $\Phi$ and $\Psi$ be strongly non-degenerate. Then the following
three conditions are equivalent:
\begin{itemize}

\item[\reflb{O-index}{\rm{(a)}}]
$\mu(\Phi^k)=\mu(\Psi^k)$ for all $k\in\N$;
\item[\reflb{O-difference}{\rm{(b)}}]
  $\mu\big(\Phi^{k+1}\big)-\mu\big(\Phi^k\big)
  =\mu\big(\Psi^{k+1}\big)-\mu\big(\Psi^k\big)$
  for all $k\in\N$;
\item[\reflb{O-signature}{\rm{(c)}}] the multiset
  $\Theta(\Phi)\sqcup\Theta(\Psi)$ breaks down into a union (perhaps,
  non-uniquely) of internally or externally matched pairs and
$$
2\lp(\Phi)+\mult_{-}(\Phi)+2\iota(\Phi) =
2\lp(\Psi)+\mult_{-}(\Psi)+2\iota(\Psi),
$$
where $\iota(\Phi)$ and $\iota(\Psi)$ stand for the number of
internally matched pairs.
\end{itemize}
\end{Theorem}

Theorem \ref{thm:index-spec} is used, for instance, in the proof of
\cite[Cor.\ F]{CGG:HZ}, which Corollary \ref{cor:ellipsoids} relies
on.

\begin{proof} The fact that \ref{O-index} is equivalent to
  \ref{O-signature} is an immediate consequence of Corollary
  \ref{cor:FM} and \eqref{eq:mu-k}.  (As in the corollary, the
  non-trivial part here is that \ref{O-index} implies
  \ref{O-signature}; the opposite direction is quite standard.)

Clearly, \ref{O-index} implies \ref{O-difference}. Assume that
\ref{O-difference} holds. Setting $\varphi=\Phi\oplus\Psi^{-1}$,
we have
  $$
  \mu\big(\varphi^{k+1}\big)-
  \mu\big(\varphi^k\big)
  =0. 
  $$
  In other words, the sequence $\mu\big(\varphi^k\big)$ is
  constant. Then we necessarily have $\mu\big(\varphi^k\big)=0$, as is
  easy to see directly or, alternatively, again by Corollary
  \ref{cor:FM} and \eqref{eq:mu-k}.  Hence, \ref{O-index} holds.
\end{proof}

\begin{Remark}
  \label{rmk:determine}
  In general, neither the index sequence $\mu\big(\Phi^k\big)$ nor the
  jump sequence $\mu\big(\Phi^{k+1}\big)-\mu\big(\Phi^k\big)$
  completely determines $\Theta(\Phi)$. For instance, rotation numbers
  forming an internally matched pair are not detected in Theorem
  \ref{thm:index-spec}.
  \end{Remark}

\section{Ultra-formal pseudo-rotations}
\label{sec:reduction}

The goal of this section and the next is to prove Theorems
\ref{thm:LA-p_leq_m} and \ref{thm:LA-ex}. The first step of the proof
is a further reduction to a combinatorial structure involving only the
mean indices of the maps $\Phi_i$ and $\Psi_j$, which we call an
ultra-formal pseudo-rotation. This is done in Theorem
\ref{cor:mean}. Then Theorems \ref{thm:numbers} and
\ref{thm:numbers-ex} are the counterparts of Theorems
\ref{thm:LA-p_leq_m} and \ref{thm:LA-ex} for ultra-formal
pseudo-rotations.

\subsection{From formal to ultra-formal}
\label{sec:UFPR}
We begin by introducing some notation. Let 
$$
\CU = \{ X_1,\ldots,X_{m+1}, Y_1,\ldots,Y_{ q } \}
$$ 
be a set of $m+ q +1$ positive real numbers. Set
\[
\rho:=\sum_iX_i+\frac{1}{2}\sum_jY_j.
\]
For $z \in\CU$, let $\CW_z$ be the multiset given by
$$
\CW_z = \big(\CU\setminus\{z\}\big)\sqcup \big\{ -Y_j/2 \, \vert \,
Y_j\neq z \big\},
$$
where the square cup stands for the union of multisets.  Thus $\CW_z$
has $m+2 q $ elements when $z=X_i$ and $m+2 q -1$ elements when
$z=Y_j$.

Throughout the rest of the paper, we call a pair of real numbers
\emph{$z$-matched} if their sum is an integer multiple of $z$.  Our
next goal is to characterize the sets formed by the reciprocals of the
mean indices of formal pseudo-rotations. To this end, it is convenient
to further formalize the concept.

\begin{Definition}
  \label{def:UFPR}
  A set $\CU=\{X_1,\ldots,X_{m+1},Y_1,\ldots,Y_{ q }\}$ of positive
  real numbers is an \emph{ultra-formal pseudo-rotation} if it meets
  the following two conditions:
  \begin{itemize}

  \item[\reflb{S-I}{\rm{(I)}}] all pairwise ratios of the elements of
    $\CU$ are irrational, except the ratios $X_i/X_i=1=Y_j/Y_j$;

  \item[\reflb{S-P}{\rm{(P)}}] for every $z\in\CU$ the multiset
    $\CW_z$ contains at least $ q $ disjoint $z$-matched pairs.
  \end{itemize}

\end{Definition}

Strictly speaking, in this definition, $\CU$ is a pair of sets
$\{X_i\}$ and $\{Y_j\}$, since the elements of these sets are treated
differently in, e.g., the definition of $\CW_z$. Furthermore, $\CW_z$
is a genuine set due to \ref{S-I}.  When $ q =0$, the matching
condition, \ref{S-P}, is void, and in this case an ultra-formal
pseudo-rotation is simply a set of $m+1$ positive numbers with
pairwise irrational ratios.

The next theorem, proved in Section \ref{sec:pf-key-lemma}, is the
first step in the proofs of Theorems \ref{thm:LA-p_leq_m} and
\ref{thm:LA-ex}, establishing a connection between the two formal
notions.

\begin{Theorem}
  \label{cor:mean}
  Let $\{\Phi_i, \Psi_j\}$ be a formal pseudo-rotation. Then the
  reciprocals $X_i=\hmu(\Phi_i)^{-1}$ and $Y_j=\hmu(\Psi_j)^{-1}$ form
  an ultra-formal pseudo-rotation such that, in addition,
 \begin{itemize}
 
 \item[\reflb{S-N}{\rm{(N)}}]
   $\sum_iX_i+\frac{1}{2}\sum_jY_j=\frac{1}{2}$.

  \end{itemize}
  Conversely, let $\CU=\{X_1,\ldots,X_{m+1},Y_1,\ldots,Y_{ q }\}$ be
  an ultra-formal pseudo-rotation satisfying \ref{S-N}. Then there
  exists a formal pseudo-rotation
  $\{\Phi_1,\ldots,\Phi_{m+1}, \Psi_1,\ldots,\Psi_{ q }\}$ with
  $\hmu(\Phi_i)=X_i^{-1}$ and $\hmu(\Psi_j)=Y_j^{-1}$.

\end{Theorem}

The reason that condition \ref{S-N} is not included in the definition
of ultra-formal pseudo-rotations is that this condition can always be
satisfied by scaling $X_i$ and $Y_j$, while \ref{S-I} and \ref{S-P}
are invariant under scaling. Thus Theorem \ref{cor:mean} asserts that
passing to the reciprocals of the mean indices gives rise to a
well-defined map from formal to ultra-formal pseudo-rotations and this
map is onto, up to rescaling.

\begin{Remark}
  For fixed $m$ and $ q $, denote by $\FF_{m, q }$ the set of formal
  pseudo-rotations and by $\CU\FF_{m, q }$ the set of ultra-formal
  pseudo-rotations up to scaling. By Theorem \ref{cor:mean}, taking
  the reciprocals of the mean indices gives rise to a well-defined map
  $\FF_{m, q }\to \CU\FF_{m, q }$, and this map is onto. However, to
  the best of our understanding, this map is rarely (almost never)
  one-to-one even if we disregard the ambiguity arising from the
  hyperbolic eigenvalues. In fact, the map fails to be one-to-one
  already when $ q =0$ and some of the ratios $\pm X_i/X_j$,
  $i\neq j$, are not distinct modulo $\Z$. This shows that on the
  formal level the non-resonance condition in Corollary
  \ref{cor:ellipsoids} cannot be removed. In other words, even when
  $ q =0$ and $m\geq 2$, there exist formal pseudo-rotations different
  from ellipsoids and hypothetically not realized by true
  pseudo-rotations.
\end{Remark}  

By \ref{cond-LSpart} and \eqref{eq:mu-hmu}, the mean indices of the
maps in Theorem \ref{thm:LA-p_leq_m} are positive. Hence, as a
consequence of Theorem \ref{cor:mean}, Theorem \ref{thm:LA-p_leq_m} is
equivalent to the following result.

\begin{Theorem}
  \label{thm:numbers}
  Let $\CU=\{X_1,\ldots,X_{m+1}, Y_1,\ldots,Y_{ q }\}$ be an
  ultra-formal pseudo-rotation. Then $ q \leq\max\{0,m-2\}$. In
  particular, $ q =0$ for $m=1,2$.
\end{Theorem}

Similarly, by Theorem \ref{cor:mean}, the next theorem is a rephrasing
of Theorem \ref{thm:LA-ex}.

\begin{Theorem}
  \label{thm:numbers-ex}
  For every integer $m\geq 3$, there exists an ultra-formal
  pseudo-rotation with $ q = 1$. Moreover, $ q $, as a function of
  $m$, is unbounded: there exist ultra-formal pseudo-rotations with
  arbitrarily large $ q $ (and hence also large $m$).
\end{Theorem}

We prove these theorems in Section \ref{sec:proofs}. 

\subsection{Proof of Theorem \ref{cor:mean}}
\label{sec:pf-key-lemma}

The proof of the theorem hinges on the following lemma.

\begin{Lemma}
  \label{lemma:key}
  Let
  $\Phi_1,\ldots,\Phi_{m+1}, \Psi_1,\ldots,\Psi_{ q } \in \tSp(2m)$ be
  strongly non-degenerate maps such that $\Phi_i$ are non-alternating
  and $\Psi_j$ are alternating. Assume that these maps have positive
  mean indices and satisfy \ref{cond-irr}. Let $\CU$ be the collection
  of the reciprocals of their mean indices: $X_i:=\hmu(\Phi_i)^{-1}$
  and $Y_j:=\hmu(\Psi_j)^{-1}$. For $\varphi =\Phi_i$ or
  $\varphi=\Psi_j$, set $z:=\hmu(\varphi)^{-1} \in\CU$ and denote by
  $\Theta(\varphi)$ the multiset formed by its rotation numbers. Then
  the conditions \ref{cond-monotone} and \ref{cond-LSpart} hold if and
  only if for every $\varphi$ the multiset
 $$
 \CW_z\sqcup\big(-z\Theta(\varphi)\big)
 $$
 breaks down into $z$-matched pairs.
\end{Lemma}

Let us now prove the theorem, assuming the lemma, which we will
establish later in this section. The main point of the argument is
eliminating the rotation numbers $\Theta(\varphi)$ from the
statement. Throughout the proofs, it will be convenient to denote by
$\Sigma(A)$ the sum of the elements of a finite multiset $A=\{a_i\}$
of real numbers: $\Sigma(A)=\sum a_i$. For instance, $\Sigma(\CW_z)$
is the sum of the elements of $\CW_z$, and we have
\begin{equation}
    \label{eq:resonance-W}
    \Sigma(\CW_z) =\begin{cases}
      \rho-z, & \text{if $z=X_i$},\\
      \rho-z/2, & \text{if $z=Y_j$}.
    \end{cases}
  \end{equation}

\begin{proof}[Proof of Theorem \ref{cor:mean}]
  First, consider a formal pseudo-rotation. By \ref{cond-LSpart} and
  \eqref{eq:mu-hmu}, the mean indices of the maps are positive and
  \ref{cond-irr} is equivalent to \ref{S-I}. To prove \ref{S-P}, fix
  $z\in\CU$ and pick $\varphi=\Phi_i$ or $\varphi=\Psi_j$ with
  $\hmu(\varphi)^{-1}=z$. By Lemma \ref{lemma:key}, the multiset
  $\CW_z\sqcup\big(-z\Theta(\varphi)\big)$ breaks down into
  $z$-matched pairs and at most $|\Theta(\varphi)|$ (the cardinality)
  of them contain an element of $-z\Theta(\varphi)$. Since
$$
|\Theta(\varphi)|=\ee(\varphi)\leq m-\hh_{-}(\varphi)
$$ 
and $\hh_{-}(\varphi)\neq 0$ when $\varphi$ is alternating, we have
$$
|\CW_z|-|\Theta(\varphi)|\geq 2 q 
$$ 
regardless of whether $z=X_i$ or $z=Y_j$. Therefore, at least $2 q $
elements of $\CW_z$ are paired among themselves. Finally, \ref{S-N} is
the resonance identity from \cite{Vi:R}. To prove \ref{S-N} directly,
it suffices to divide \eqref{eq:mu-count} below (see Step 1 in the
proof of Lemma \ref{lemma:key}) by $k$, let $k\to\infty$ and use the
facts that $\CN(t)/t\to \rho$, where $\CN(t)$ is the counting function
defined by \eqref{eq:counting}, and
$\mu(\varphi^k)/k\to\hmu(\varphi)$. Alternatively, \ref{S-N} follows
from \ref{cond-LSpart} alone; see \cite{Vi:R} or \cite{GiK}.

Conversely, suppose that $\CU$ is an ultra-formal pseudo-rotation
satisfying \ref{S-N}. Fix $z\in\CU$ and $ q $ disjoint $z$-matched
pairs in $\CW_z$. There are $|\CW_z|-2 q $ remaining elements
$\alpha_{r}$, $r=1,\ldots,|\CW_z|-2 q $. More precisely, we have $m$
such elements if $z=X_i$ and $m-1$ if $z=Y_j$. Let $T_z$ be the
multiset of the fractional parts $\{\alpha_{r}/z\}$ of the ratios
$\alpha_{r}/z$.  (These are irrational by \ref{S-I}.) Let $\varphi_z$
be the direct sum of the rotations $R_\theta$ with $\theta\in T_z$. If
$z=Y_j$, we add one negative hyperbolic block to $\varphi_z$. Then
$\varphi_z$ is strongly non-degenerate with $\Theta(\varphi_z)=T_z$,
non-alternating if $z=X_i$ and alternating if $z=Y_j$. It remains to
compose $\varphi_z$ with a loop so that $\hmu(\varphi_z)=z^{-1}$. By
\eqref{eq:mean}, this amounts to checking that $z^{-1}-2\Sigma(T_z)$
is an even integer if $z=X_i$ and an odd one if $z=Y_j$.

The $z$-matched pairs given by \ref{S-P} together with
$\big(\alpha_{r},-z\{\alpha_{r}/z\}\big)$ break down
$\CW_z\sqcup(-zT_z)$ into disjoint $z$-matched pairs. Summing them all
up and dividing by $z$, we see that
\begin{equation}
\label{eq:mean-cor-pf}
z^{-1}\Sigma(\CW_z)-\Sigma(T_z) \in \Z. 
\end{equation}
Next, from \eqref{eq:resonance-W}, \ref{S-N} and
\eqref{eq:mean-cor-pf} we conclude that
 $$
 z^{-1}/2-\Sigma(T_z) \in \Z\ \ \text{if $z=X_i$}\quad\text{and}\quad
 z^{-1}/2-1/2-\Sigma(T_z) \in \Z\ \ \text{if $z=Y_j$}.
 $$
 Hence, $z^{-1}-2\Sigma(T_z)$ is even if $z=X_i$ and odd if $z=Y_j$.
 Finally, \ref{cond-irr} is \ref{S-I} and, by Lemma \ref{lemma:key},
 the maps $\varphi_z$ satisfy \ref{cond-monotone} and
 \ref{cond-LSpart}. As $z$ varies through $\CU$, we obtain the desired
 formal pseudo-rotation.
\end{proof}

\begin{Remark}
  It is useful to keep in mind that the maps $\varphi_z$ from the
  proof of Theorem \ref{cor:mean} are semisimple, and they are
  elliptic for $z=X_i$ and have $m-1$ elliptic blocks and one negative
  hyperbolic block for $z=Y_j$.
\end{Remark}

\begin{proof}[Proof of Lemma \ref{lemma:key}]
  The idea of the proof is that we can determine the position of an
  iterate $\varphi^k$ on the list of all iterates, ordered by the mean
  index, in two different ways.  The first way is by the
  Conley--Zehnder index and the other is via the number of good
  iterates with smaller or equal mean index. The latter is the
  counting function $\CN(t)$ given by \eqref{eq:counting}. By
  \ref{cond-monotone} and \ref{cond-LSpart}, the two ways are
  equivalent, resulting in a relation, \eqref{eq:mu-count}, between
  $\CN$ and the index.  This is Step 1 of the proof. In Step 2, we
  rephrase this relation as an identity on certain generalized
  Beatty-type sequences in the framework of Corollary \ref{cor:FM} and
  frequency matching; see \eqref{eq:final}. The frequencies are
  essentially the ratios of the mean indices and the rotation
  numbers. Applying the corollary to this identity in Step 3, we
  obtain the desired matched pairs.
 
  \smallskip\noindent\emph{Step 1:} For $t>0$, define the counting
  function $\CN(t)$ by
\begin{equation}
  \label{eq:counting}
  \CN(t):=\sum_{i}\Big\lfloor \frac{t}{\hmu(\Phi_i)}\Big\rfloor
  +\sum_{j}\Big(\Big\lfloor \frac{t}{\hmu(\Psi_j)}\Big\rfloor
  -\Big\lfloor \frac{t}{2\hmu(\Psi_j)}\Big\rfloor\Big).
\end{equation}
The function $\CN(t)$ counts the number of good iterates of $\Phi_i$
and $\Psi_j$ with mean index at most $t$: the iterates $\Phi_i^k$ with
$k\hmu(\Phi_i)\leq t$ and the odd iterates $\Psi_j^k$ with
$k\hmu(\Psi_j)\leq t$. We claim that the conditions
\ref{cond-monotone} and \ref{cond-LSpart} hold if and only if for
every $\varphi\in\{\Phi_i,\Psi_j\}$ and $k\in\N$, we have
\begin{equation}
    \label{eq:mu-count}
    \mu\big(\varphi^k\big)  =\begin{cases}
      m+2\CN\big(\hmu(\varphi^k)\big),
      & \text{if $\varphi^k$ is good},\\
      m+2\CN\big(\hmu(\varphi^k)\big) +1,
      & \text{if $\varphi^k$ is bad}.
    \end{cases}
\end{equation}
Namely, first suppose that \ref{cond-monotone} and \ref{cond-LSpart}
hold. Then, for good iterates, the index strictly increases with mean
index and fills $m+2,m+4,\ldots$. This explains the first line. As for
the second, by \ref{cond-irr}, the mean index of a bad iterate is
different from the mean indices of all good iterates. Hence, condition
\ref{cond-monotone} implies that the index of a bad iterate is
squeezed between the indices of the nearest good iterates on both
sides. (Recall that $\Psi_j^1$ is a good iterate.) The squeeze is
strict for parity reasons: the index of an even iterate of an
alternating map has parity $m+1$. Therefore, the index is one more
than $m+2\CN\big(\hmu(\varphi^k)\big)$. This is the second line of
\eqref{eq:mu-count}.

For the converse, suppose that \eqref{eq:mu-count} holds. By
\ref{cond-irr}, the counting function $\CN(t)$ jumps exactly by one at
the mean indices of good iterates, and does not change at the bad
ones. This, together with \eqref{eq:mu-count}, implies
\ref{cond-LSpart}. As for monotonicity \ref{cond-monotone}, let
$\varphi$ and $\varphi'$ be two iterates with
$\hmu(\varphi) < \hmu(\varphi')$. If one of them is good, or if the
interval in between contains a good iterate, we have
$\mu(\varphi) < \mu(\varphi')$ by \eqref{eq:mu-count}. Otherwise,
$\mu(\varphi) = \mu(\varphi')$ by the discussion above and
\eqref{eq:mu-count} again. To summarize, conditions
\ref{cond-monotone} and \ref{cond-LSpart} hold if and only if
\eqref{eq:mu-count} is satisfied.

\smallskip\noindent\emph{Step 2:} Fix $\varphi\in\{\Phi_i,\Psi_j\}$. By
\eqref{eq:mean} and \eqref{eq:mu-k}, we have
\begin{equation}
  \label{eq:mu-k-phi}
  \mu\big(\varphi^k\big)=\ee(\varphi)
  +\big(\hmu(\varphi)-2\Sigma\big(\Theta(\varphi)\big)\big)k
  +2b_{\Theta(\varphi)}(k)
\end{equation}
where $b_{\Theta(\varphi)}(k)=\sum_l\lfloor k\theta_l(\varphi)\rfloor$
as in Corollary \ref{cor:FM}. Our goal is to set \eqref{eq:mu-k-phi}
equal to \eqref{eq:mu-count} with $t=k\hmu(\varphi)$ and then apply
frequency matching. However, to do so, we need to modify
\eqref{eq:mu-count} slightly. Namely, as it is, the counting function
contains Beatty-type terms with a negative sign (from the even
iterates of the $\Psi_j$) and the term involving $\varphi$ itself in
$\CN\big(k\hmu(\varphi)\big)$ has rational frequency. Hence, we first
separate the latter term and then collect all Beatty-type terms on one
side.

The term of $\varphi$ in $\CN\big(k\hmu(\varphi)\big)$ is $k$ if
$\varphi$ is non-alternating and $k-\lfloor k/2\rfloor$ if $\varphi$
is alternating. In the latter case, $2(k-\lfloor k/2\rfloor)$ plus $1$
for even $k$ (the bad iterates) in \eqref{eq:mu-count} gives $k+1$ for
every $k$. Therefore, the contribution of $\varphi$ to the right-hand
side of \eqref{eq:mu-count} is $2k$, and respectively $k+1$, in the
non-alternating/alternating cases. The terms coming from the maps
other than $\varphi$ add up to
$2\big(b_{A(\varphi)}(k)-b_{B(\varphi)}(k)\big)$, where
\[
A(\varphi)
:=\Big\{\frac{\hmu(\varphi)}{\hmu(\varphi')}\ \Big\vert\
\varphi'\in\{\Phi_i,\Psi_j\},\ \varphi'\neq\varphi\Big\},
\quad
B(\varphi)
:=\Big\{\frac{\hmu(\varphi)}{2\hmu(\Psi_j)}\ \Big\vert\
\Psi_j\neq\varphi\Big\}.
\]
The elements of $A(\varphi)$, $B(\varphi)$ and $\Theta(\varphi)$ are
positive and irrational. For $A(\varphi)$ and $B(\varphi)$, this
follows from \ref{cond-irr} and the condition that the mean indices
are positive. For $\Theta(\varphi)$, this is a consequence of strong
non-degeneracy and the definition of the rotation numbers. We conclude
that \eqref{eq:mu-count} is equivalent to
\begin{equation}
    \label{eq:enter-beatty}
    \mu\big(\varphi^k\big)  =\begin{cases}
      m+2k+2b_{A(\varphi)}(k)-2b_{B(\varphi)}(k),
      & \text{if $\varphi$ is non-alternating},\\
      m+1+k+2b_{A(\varphi)}(k)-2b_{B(\varphi)}(k),
      & \text{if $\varphi$ is alternating}
    \end{cases}
  \end{equation}
for all $k\in\N$.

Next, we combine \eqref{eq:mu-k-phi} and \eqref{eq:enter-beatty}, and
collect the Beatty-type sequences on one side. We conclude that
\eqref{eq:mu-count} is equivalent to
\begin{equation}
  \label{eq:final}
 2\big(b_{A(\varphi)}-b_{B(\varphi)}-b_{\Theta(\varphi)}\big)(k)=
  \begin{cases}
    \ee(\varphi)-m
    +\big(\hmu(\varphi)-2\Sigma(\Theta(\varphi))-2\big)k,
    \\
    \ee(\varphi)-m-1
    +\big(\hmu(\varphi)-2\Sigma(\Theta(\varphi))-1\big)k
  \end{cases}
\end{equation}
for all $k\in\N$, where the first line is for non-alternating
$\varphi$ and the second one is for alternating $\varphi$.

\smallskip\noindent\emph{Step 3:} With \eqref{eq:final} in place, we
are in the setting of Corollary \ref{cor:FM} with $A(\varphi)$ and
$B(\varphi)\sqcup\Theta(\varphi)$ playing the roles of $A$ and
$B$. Set $z:=\hmu(\varphi)^{-1}$ and observe that
\begin{equation}
\label{eq:W-AB}
\hmu(\varphi)\CW_{z}=A(\varphi)\sqcup\big(-B(\varphi)\big).
\end{equation} 
Hence $A(\varphi)\sqcup B(\varphi)\sqcup\Theta(\varphi)$ breaks down
into internally or externally matched pairs if and only if
$$
\CW_{z}\sqcup\big(-z\,\Theta(\varphi)\big)
$$ 
breaks down into $z$-matched pairs. Indeed, by \eqref{eq:W-AB},
multiplying this multiset by $\hmu(\varphi)=z^{-1}$ turns it into
$A(\varphi)\sqcup(-B(\varphi))\sqcup(-\Theta(\varphi))$, and a pair in
that union has integer sum exactly when the corresponding pair in
$A(\varphi)\sqcup B(\varphi)\sqcup\Theta(\varphi)$ is internally
matched (i.e., both elements are in $A(\varphi)$ or in
$B(\varphi)\sqcup\Theta(\varphi)$) or externally matched (i.e., one
element is in each). As a consequence, and by Steps 1--2, it suffices
to show that \eqref{eq:final} holds for all $\varphi$ if and only if
the former breakdown exists for all $\varphi$.

If \eqref{eq:final} holds, then the LHS of \eqref{eq:final} is a
linear function of $k$, and the matching follows from Corollary
\ref{cor:FM}.  Conversely, suppose that for every $\varphi$ the
multiset $A(\varphi)\sqcup B(\varphi)\sqcup\Theta(\varphi)$ breaks
down into internally or externally matched pairs. By Corollary
\ref{cor:FM}, we have
$$
b_{A(\varphi)}(k)-b_{B(\varphi)}(k)-b_{\Theta(\varphi)}(k)=ak-c
$$
for all $k\in\N$. Furthermore, the coefficients $a$ and $c$ are
integers, and
\[
  a=\Sigma(A(\varphi))-\Sigma(B(\varphi))-
  \Sigma(\Theta(\varphi))\quad\text{and}
  \quad
c=\big(|A(\varphi)|-|B(\varphi)|-|\Theta(\varphi)|\big)/2.
\]
We need to show that $2(ak-c)$ is the right-hand side of
\eqref{eq:final}. For the constant term, we note that
$$
|A(\varphi)|=m+ q , \quad |\Theta(\varphi)|
=\ee(\varphi), \quad
|B(\varphi)|= q  \text{ or }  q -1
$$
when $\varphi$ is non-alternating/alternating. Thus, $-2c$ is
$\ee(\varphi)-m$ and, respectively, $\ee(\varphi)-m-1$, as
required. For the linear term, using \eqref{eq:W-AB} and
\eqref{eq:resonance-W}, we write
\begin{align*}
  \hmu(\varphi)  -2\Sigma(A(\varphi)) + 2\Sigma(B(\varphi))
  &= \hmu(\varphi)-2 \hmu(\varphi) \Sigma(\CW_{z}) \\ 
  &= 2 + \hmu(\varphi)( 1- 2\rho)  \text{   or  } 1
    + \hmu(\varphi)( 1- 2\rho),
\end{align*}
depending on whether $\varphi$ is non-alternating or alternating. It
remains to show that the second term on the RHS, i.e.,
$\hmu(\varphi)( 1- 2\rho)$, vanishes for all $\varphi$. In other
words, we need to prove that $\rho=1/2$. To this end, first observe
that the RHS is an integer for all $\varphi$. This follows from the
fact that the LHS is the difference of two integers:
$$
\big(\hmu(\varphi)-2\Sigma(\Theta(\varphi))\big)-2a;
$$ 
see \eqref{eq:mean}. We conclude that 
$$
\hmu(\varphi)( 1- 2\rho) \in \Z
$$ 
for all $\varphi$.  By \ref{cond-irr}, this is only possible if
$1- 2\rho=0$. This completes the proof of Lemma \ref{lemma:key} and
hence of Theorem \ref{cor:mean}.
\end{proof}

\section{Proofs of Theorems \ref{thm:numbers} and
  \ref{thm:numbers-ex}}
\label{sec:proofs}

In this section, we prove Theorems \ref{thm:numbers} and
\ref{thm:numbers-ex}. Throughout the section, $\CU$ is an ultra-formal
pseudo-rotation.

\subsection{Preparation}
\label{sec:prep}
For $x\in\CU$, we call an $x$-matched pair $a, b \in \CU$ such that
$a+b=x$ an \emph{exact decomposition} of $x$, and refer to $a$ and $b$
as the \emph{summands} of $x$. Clearly, then $a<x$ and $b<x$, and in
any collection of distinct exact decompositions of the same element
$x$ all summands are distinct.

Here is a rough outline of the proof. For a suitably chosen element
$w \in \CU$, every $w$-matched pair $a, b \in \CW_w$ is an exact
decomposition. By \ref{S-P}, there are at least $ q $ disjoint
$w$-matched pairs, and hence at least $2 q $ summands of $w$ in
$\CU\setminus\{ w \}$. In particular, the number of
\emph{non-summands} of $w$ in $\CU \setminus \{ w \}$ is at most
$m- q $.  To finish the proof, it now suffices to show that
$m- q \ge 2$.  To this end, consider the vector space
$V=\mathrm{span}_\Q\,\CU \subset \R$ and the $\Q$-linear maps
$g \colon V \to \R$. We note that $\dim V \ge 2 $ by \ref{S-I} and use
this fact to construct a map $g$ such that a certain element $w$
associated with $g$, which we call the extremizer, has at least two
non-summands in $\CU \setminus \{ w \}$. (This $\Q$-linear map is
either $\id$ or $\beta(x) = f(x)-s(S_{\min})\,x$ or the map
$s(S_{\max})x-f(x)$, which is constructed similarly to $\beta $ but
with $f$ replaced by $-f$.)

To illustrate the idea, let us prove the weaker upper bound
$ q \leq m$. Note that for $z=\max \CU$ all $z$-matched pairs are
exact decompositions. The reason is that $0<a+b<2z$ when
$a, b\in \CW_z$ are positive, $-z<a+b<z$ when one of them is negative
(or both), and \ref{S-I} rules out $a+b=0$. Since $\CU$ is an
ultra-formal pseudo-rotation, there are at least $ q $ such pairs and
hence at least $ q $ exact decompositions. Their summands are distinct
elements of $\CU$ different from $z$. Thus
$m+ q =|\CU\setminus\{z\}|\geq 2 q $ and hence $m\geq q $. This
argument also suggests that the property of having only exact
decompositions as matched pairs is associated with extremizing
features of $z$. This is formalized by the following lemma.

Let $V\subset\R$ be the $\Q$-span of $\CU$.  Note that
$1<\dim V \leq m+ q +1$ by \ref{S-I}, where $\dim V$ is the dimension
of $V$ over $\Q$.

\begin{Lemma}
  \label{lemma:max}
  Let $g\colon V\to\R$ be a $\Q$-linear map such that $g\geq 0$ on
  $\CU$ with $ g \not\equiv 0$ and $g(Y_j)>0$ for all $j$. Let $w$ be
  the smallest element of $\CU$ such that $g(w)=\max_{\CU} g$. Then
  every $w$-matched pair in $\CW_w$ is an exact decomposition of $w$
  and $g>0$ on both summands. In particular, $w$ has at least $2 q $
  distinct summands in $\CU\setminus\{w\}$ with $g>0$.
\end{Lemma}

We call $w$ the \emph{extremizer} for $g$ and write $w=w_g$ when $w$
meets the conditions of the lemma: $w\in\CU$ is the smallest element
with $g(w)=\max_{\CU} g$. For instance, $z=\max \CU$ is the extremizer
for $g=\id$.

\begin{proof}
  Assume that $a+b=kw$ for some $k\in\Z$ and $a,b\in\CW_w$. Let us
  first show that $a$ and $b$ have the same sign. To this end,
  applying $g$, we obtain
\begin{equation}
\label{eq:lem-max}
g(a)+g(b)=kg(w)=k\max_{\CU} g =: kM
\end{equation}
and also note that
$$
0\leq g\leq M \textrm{ on $\CU$ and }
-M/2\leq g(-Y_j/2)<0 \textrm{ for all $j$}.
$$
These two facts imply that $a$ and $b$ have the same sign. Indeed, if
the signs were opposite, we would have $k=0$. This is impossible by
\ref{S-I}, since then $\CU$ would have an element equal to $Y_j/2$.

Furthermore, $a>0$ and $b>0$. To see this, assume the contrary, i.e.,
$a=-Y_i/2$ and $b=-Y_j/2$ for some $i\neq j$. Thus $k=-1$ by the
second inequality above. Then, \eqref{eq:lem-max} implies that
$g(-2a)=g(-2b)=M$, and hence $Y_i=-2a>w$ and $Y_j= -2b >w$ due to the
choice of $w$ as the smallest element on which $g$ attains its maximum
value $M$.  This is incompatible with the condition that $a+b=-w$.

Thus $a$ and $b$ are both positive, i.e., $a,b \in \CU$. Then either
$k=1$ or $k=2$. Arguing as above, we see that, in the latter case, we
must have $g(a)=g(b) =M$. Hence $a>w$ and $b >w$, which is impossible
since $a+b=2w$. We conclude that $a,b \in \CU$ and $a+b =w$.

To show that $g(a)>0$ and $g(b) > 0$, recall first that $g\geq 0$ on
$\CU$. Thus the only alternative is that $g(a)=0$ and/or $g(b)=0$. If,
say, $g(a)=0$, we must have $g(b)=M$, and hence $b>w$, which is again
impossible since $a+b=w$ and both $a$ and $b$ are positive.
\end{proof}

Next, let $f\colon V\to\R$ be a $\Q$-linear map such that the slopes
$$
s(x):=f(x)/x, \quad x\in\CU
$$
are pairwise distinct. The existence of such a map follows from
\ref{S-I}. Namely, fix a pair $x\neq y$ in $\CU$. The set of
$\Q$-linear maps $V\to\R$ is a real vector space of dimension equal to
$\dim V > 1$, and the ones giving the same slope at $x$ and $y$ form a
linear subspace since the condition $y f(x)=x f(y)$ is linear in
$f$. This subspace has positive codimension. Indeed, for a linear map
$f$ with $f(x)=0$ and $f(y)=1$, which exists by \ref{S-I}, the slopes
are distinct: $s(x)=0 \neq s(y)=1/y$. Since there are finitely many
pairs, one can always choose a function $f$ in the complement of these
subspaces. Next, we fix such a map $f$ and denote by $S_{\min}$ and
$S_{\max}$ the elements in $\CU$ with the smallest and the largest
slopes.  Define a $\Q$-linear function $\beta\colon V\to\R$ by
\begin{equation}
  \label{eq:beta}
\beta(x):=f(x)-s(S_{\min})x;
\end{equation}
see Figure \ref{Fig}. Note that for $x\in\CU$, we have
$\beta(x) = x\big(s(x)-s(S_{\min})\big)$. It follows that
$\beta\geq 0$ on $\CU$, with equality only at $S_{\min}$. As a
consequence, $\beta \neq 0 $ on $\CW_{S_{\min}}$ and has the same sign
as its argument, i.e., $\beta(x)$ and $x$ have the same sign for
$x\in\CW_{S_{\min}}$. In particular, $\beta(x)<0$ for
$x\in\CW_{S_{\min}}$ if and only if $x$ is one of the elements
$-Y_j/2$.

\begin{Lemma}
  \label{lemma:slope}
  Both $S_{\min}$ and $S_{\max}$ are among the $X_i$.
\end{Lemma}

\begin{proof}
  Without loss of generality, we may assume that $ q \geq 1$; for
  otherwise the statement is true automatically.  Let us first focus
  on $S_{\min}$. By \ref{S-P}, $\CW_{S_{\min}}$ contains $ q $
  disjoint $S_{\min}$-matched pairs. We claim that each pair contains
  exactly one of the elements $-Y_j/2$. This would complete the proof
  since having $S_{\min}=Y_j$ implies that there are at most $ q -1$
  such elements in $\CW_{S_{\min}}$, while we need $ q $ disjoint
  matched pairs. Next, let us prove the claim. Consider an
  $S_{\min}$-matched pair, i.e., $x,y\in\CW_{S_{\min}}$ with
  $x+y=kS_{\min}$ for some $k\in\Z$. Applying $\beta$ and using the
  fact that $\beta(S_{\min})=0$, we obtain
$$
\beta(x)+\beta(y)=0.
$$
Then, as we have pointed out above, $x$ and $y$ have different signs,
and hence exactly one of them is of the form
$-Y_j/2 \in \CW_{S_{\min}}$. For $S_{\max}$, we use the same argument
but with $f$ replaced by $-f$.
\end{proof}

An important consequence of Lemma \ref{lemma:slope} is that $g=\beta$
satisfies the requirements of Lemma \ref{lemma:max}. Indeed,
$\beta\geq 0$ by construction and the condition that $\beta(Y_j)>0$
holds by Lemma \ref{lemma:slope} since $\beta > 0$ on
$\CU\setminus \{ S_{\min} \}$.

\begin{Remark}
  The proof of Lemma \ref{lemma:slope} has an interesting application
  to dynamics. Namely, as a consequence of the lemma and its proof,
  there are no non-degenerate Reeb pseudo-rotations with $ q \neq 0$
  that are sufficiently $C^1$-close to an irrational ellipsoid, i.e.,
  such that the Reeb vector field is $C^1$-close to the one on an
  irrational ellipsoid. (Hypothetically, such a pseudo-rotation could
  be close to a rational one.) Indeed, assume that such a
  pseudo-rotation exists. Then the mean indices of the alternating
  orbits must be large, and hence their reciprocals $Y_j$ are
  small. By the proof of the lemma, we have an $S_{\min}$-matched pair
  of the form $X_i-Y_j/2=kS_{\min}$, where $S_{\min}=X_{i'}$. This is
  impossible when $Y_j<\eps$, where $\eps=\min |X_i- kX_{i'}|>0$ with
  the minimum taken over all $i \neq i'$ and $k$.
\end{Remark}  

\subsection{Proof of Theorem \ref{thm:numbers}}
\label{sec:pf-numbers}
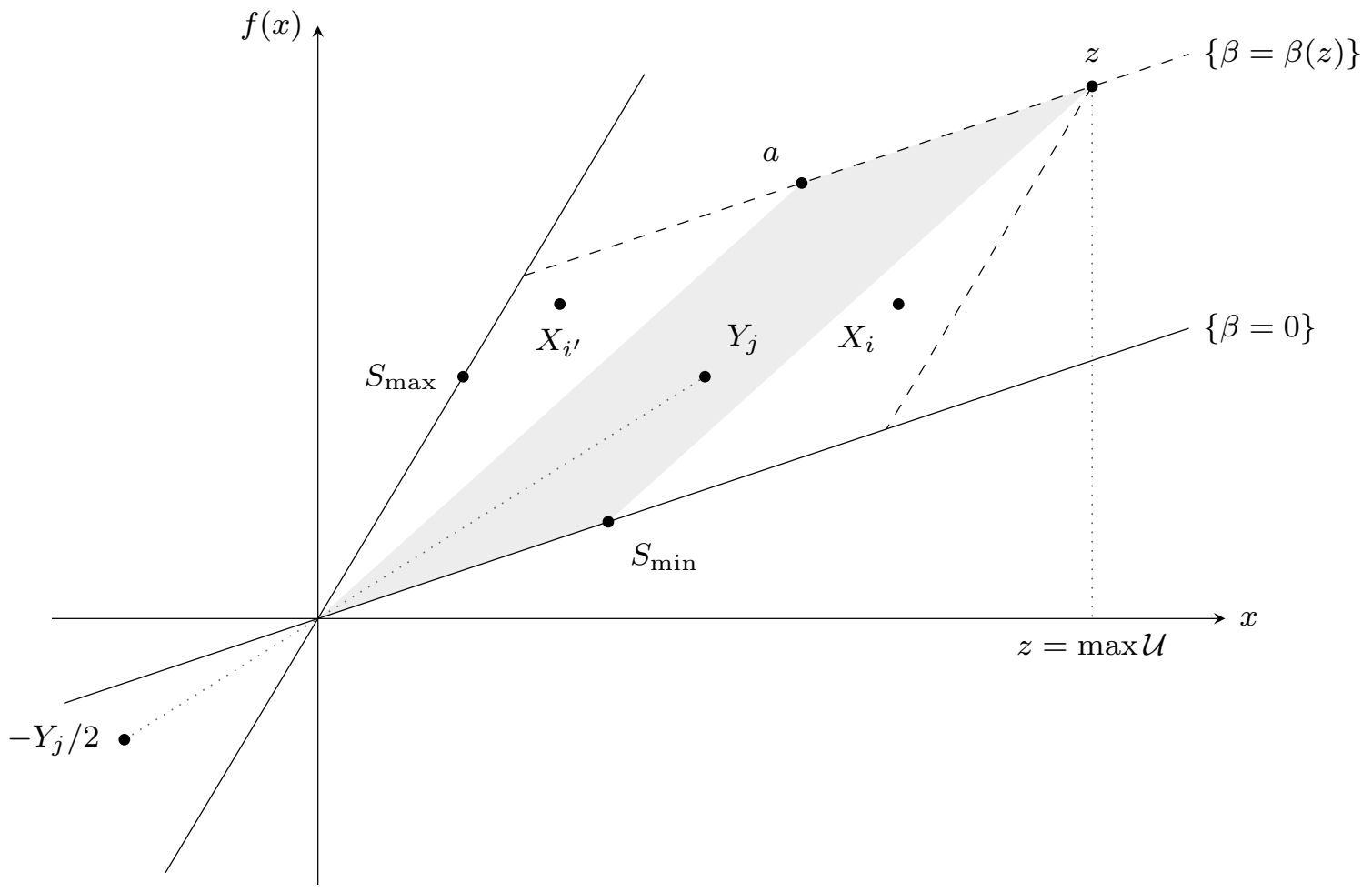
\begin{figure}[htbp]
  \centering
  \scalebox{0.8}{%
    \begin{tikzpicture}[x=2.2cm,y=2.2cm,line width=0.3pt,>=stealth,
        every node/.style={font=\scriptsize},
        pt/.style={circle,fill=black,inner sep=0pt,minimum size=3pt}]
      \fill[black!7] (0,0) -- (1.2,0.4) -- (3.2,2.2) -- (2.0,1.8) -- cycle;
      \draw[->] (-1.1,0) -- (3.75,0) node[right] {$x$};
      \draw[->] (0,-1.1) -- (0,2.45) node[left] {$f(x)$};
      \draw (-1.05,-0.35) -- (3.6,1.2) node[right] {$\{\beta=0\}$};
      \draw (-0.63,-1.05) -- (1.35,2.25);
      \draw[dashed] (0.85,1.417) -- (3.6,2.333) node[right] {$\{\beta=\beta(z)\}$};
      \draw[dotted,line width=0.4pt,black!60] (3.2,2.2) -- (3.2,0) node[below,black] {$z=\max\CU$};
      \draw[dashed] (3.2,2.2) -- (2.35,0.783);
      \draw[dotted,line width=0.4pt,black!60] (-0.8,-0.5) -- (1.6,1.0);
      \node[pt,label={[label distance=1pt]below right:$S_{\min}$}] at (1.2,0.4) {};
      \node[pt,label={[label distance=1pt]left:$S_{\max}$}] at (0.6,1.0) {};
      \node[pt,label={[label distance=1pt]above left:$a$}] at (2.0,1.8) {};
      \node[pt,label={[label distance=1pt]above:$z$}] at (3.2,2.2) {};
      \node[pt,label={[label distance=1pt]below left:$X_i$}] at (2.4,1.3) {};
      \node[pt,label={[label distance=1pt]below:$X_{i'}$}] at (1.0,1.3) {};
      \node[pt,label={[label distance=1pt]above right:$Y_j$}] at (1.6,1.0) {};
      \node[pt,label={[label distance=1pt]left:$-Y_j/2$}] at (-0.8,-0.5) {};
    \end{tikzpicture}%
  }
  \caption{The figure shows the points $(x,f(x))$, $x\in\CU$, which
    are used to define the slope function $s(x)$: $s(x)$ for
    $x\in \CU$ is the slope of the ray passing through $(x, f(x))$.}
  \label{Fig}
\end{figure}

Arguing by contradiction, we assume throughout the proof that
$ q \geq 1$ and $ q \geq m-1$ and use this condition
repeatedly. Recall that $g=\beta$ meets the requirements of Lemma
\ref{lemma:max}. Let $w=w_\beta$ be its extremizer. We will show in
Step 3 that at least two elements in $\CU\setminus\{w\}$ are not
summands of $w$. Hence, $w$ has at most $m+ q -2\leq 2 q -1$
summands. On the other hand, $w$ must have at least $2 q $ summands by
Lemma \ref{lemma:max} (see Step 1), leading to a contradiction.
Figure \ref{Fig} visualizes some points of the proof.

\smallskip\noindent\emph{Step 1: Let $w$ be an extremizer as in Lemma
  \ref{lemma:max}. Then at most one element in $\CU\setminus\{w\}$ is
  not a summand of $w$.}  Indeed, by our assumption on $ q $, we have
  $$
|\CU\setminus\{w\}| = m+ q  \leq 2 q +1,
  $$
  and, by Lemma \ref{lemma:max}, at least $2 q $ elements in this set
  are distinct summands of $w$.

  \smallskip \noindent\emph{Step 2: At least one of the elements
    $S_{\min}$ and $S_{\max}$ is a summand of $z=\max \CU$.}  To prove
  this, by Step 1, it suffices to show that both $S_{\min}$ and
  $S_{\max}$ are distinct from $z$ since $z$ is the extremizer for
  $g=\id$.  This, in turn, is a consequence of the fact that neither
  $S_{\min}$ nor $S_{\max}$ admits an exact decomposition. Indeed,
  suppose that $a+b=S_{\min}$ is an exact decomposition of
  $S_{\min}$. Then
  $$
  \beta(a)+\beta(b)=\beta(S_{\min})=0
  $$
  by the definition of $\beta$ from \eqref{eq:beta}. This is
  impossible since $\beta>0$ on $\CU\setminus\{ S_{\min}\}$ and hence
  $\beta(a)+\beta(b)>0$. Thus $z \neq S_{\min}$. Replacing $f$ by
  $-f$, we see that $z \neq S_{\max}$. Therefore, at least one of
  $S_{\min}$ and $S_{\max}$ must be a summand.

  \smallskip\noindent\emph{Step 3: The maximum $z=\max\CU$ and at
    least one of the elements $S_{\min}$ and $S_{\max}$ are not
    summands of $w=w_\beta$ and are different from $w$; this
    contradicts Step 1.}  Without loss of generality, we may assume
  that $S_{\min}$ is a summand of $z$. (Otherwise, we replace $f$ with
  $-f$.)  Then
  $$
  a +S_{\min} = z
  $$
  for some $a \in \CU$. In particular, $a<z$.  Applying $\beta$ to
  this decomposition, we see that $\beta(a)=\beta(z)$. We cannot have
  $z=w$; for then we would also have $\beta(a)=\beta(z)=\beta(w)$ and
  $a<z=w$, contradicting the choice of $w$ as the smallest element on
  which $\beta$ attains its maximum. Therefore, $w<z=\max \CU$.
  Hence, $z$ cannot be a summand of $w$.
  
  Clearly, $S_{\min}\neq w$ since $\beta(S_{\min})=0$, whereas
  $\beta(w)>0$.  To see that $S_{\min}$ cannot be a summand of $w$,
  assume the contrary: $b+S_{\min}=w$ for some $b \in \CU$ and, in
  particular, $b<w$. Again, applying $\beta$ to this identity, we have
  $\beta(b)=\beta(w)$, which is impossible since $w$ is an extremizer
  for $\beta$.  This completes the proof of the theorem.\qed

\subsection{Proof of Theorem \ref{thm:numbers-ex}}
\label{sec:pf-numbers-ex}
The proof of the theorem amounts to constructing specific examples of
ultra-formal pseudo-rotations with $ q =1$ or $ q $ arbitrarily large.

\subsubsection{Examples with $m\geq 3$ and $ q =1$.}
\label{sec:p=1}
Fix $m\geq 3$ and $0<a<b$ with $a/b$ irrational. Set
$$
  X_i:=a+(i-1)b, \quad X_{m+1}:=b-a, \quad Y:=2b
$$
where $1\leq i\leq m$. We claim that this collection satisfies
\ref{S-I} and \ref{S-P}. The former follows from the fact that the
ratio $a/b$ is irrational. To prove \ref{S-P}, we need to find one
$z$-matched pair in $\CW_z$ for each $z\in\CU$. These are
\begin{align*}
  X_{i+1}-\tfrac{1}{2}Y=X_i,\quad X_{m-2}+Y
  = X_m,\quad X_1-\tfrac{1}{2}Y=-X_{m+1}, \quad X_2+X_{m+1}=Y
\end{align*}
where $1\leq i\leq m-1$ in the first identity. Note that we used
$m\geq 3$ in the second identity.

\subsubsection{Examples with arbitrarily large $ q $}
It remains to find an ultra-formal pseudo-rotation $\{X_i, Y_j\}$ with
arbitrarily large $ q = | \{ Y_j \} |$.  Fix $ q > 0$ and choose
$Y_1, \dots, Y_{ q }$ so that
$$
1, Y_1, \dots, Y_{ q }
$$ 
are linearly independent over $\Q$ and $Y_j >2$ for all $j$. Consider
the direct product $B$ of $ q $ copies of the set
$\{0,1,\ldots, q \}$:
$$
B := \{0,1,\ldots,  q \}^{ q }.
$$
Denote by $e_j$ the $j$'th ``coordinate vector'' in $B$, i.e.,
$$
e_j = (0, \dots, 1, \dots, 0) \in B
$$ 
where $1$ is in the $j$'th slot. Define $( q +1)^{ q } $ real numbers
$A_b \in \R$ indexed by $b\in B$ by requiring that $A_0 =0$ and
$$
A_{b +e_j} = A_b + Y_j/2. 
$$
Finally, set $X^-_b = A_b -1$ and $X^+_b = A_b +1$. Note that
$X^\pm_b>0$ for $b \in B\setminus \{0\}$ since $Y_j>2$.

We claim that for $ q >2$ the collection $Y_j$ and $X_b^{\pm}$ for
$b \in B\setminus \{0\}$ is an ultra-formal pseudo-rotation.  (We did
not include $X^{\pm}_0$ to satisfy the positivity requirement.)
Condition \ref{S-I} readily follows from the assumption that
$1,Y_1,\ldots, Y_{ q }$ are linearly independent over $\Q$. Let us
prove \ref{S-P}.

We have
$$
X^+_b + X^-_b = 2A_b
$$
for all $b \in B$. Observe that, for
$b \in \{e_j, 2e_j, \dots, q e_j\}$, the RHS is an integer multiple of
$Y_j$. More precisely, $2A_{re_j}=rY_j$. Hence this gives $ q $
disjoint $Y_j$-matched pairs in $\CW_{Y_j}$ for each $j$.

Next, let us fix $X^\pm_b$. Assume first that $b+e_j \in B$ for all
$j$. Then we have the following $ q $ disjoint $X^\pm_b$-matched
pairs:
$$
X^\pm_{b+e_j} -  Y_j/2 = X^\pm_b. 
$$
If $b+e_j \notin B$ for some $j$, then $b-2e_j \in B\setminus\{0\}$
due to the assumption that $ q >2$. (This is the only point of the
argument where this requirement is used.) In this case, we use instead
the following pair
$$
X^\pm_{b-2e_j} + Y_j = X^\pm_b. 
$$
This proves \ref{S-P} for $ q >2$.

The case of $ q =1$ is covered in Section \ref{sec:p=1}. Furthermore,
for both $ q =1$ and $ q =2$ we can take $B=\{ 0, 1, 2, 3\}^{ q }$ as
the index set. Then the same construction gives an example with
$ m+1=\vert \{X_i\} \vert = 2(4^{ q } -1)$. (The only role of the
$ q >2$ assumption was in the backward shift of the last step of the
argument.)  This completes the proof of Theorem
\ref{thm:numbers-ex}. \qed

 \begin{Remark}
The above construction requires
$$
m+1=\vert \{X_i\} \vert= 2\big(( q +1)^{ q }-1\big )
$$
when $ q \geq 3$. We believe that this is far from optimal, but we do
not know the minimal $m$ as a function of $ q $.

\end{Remark}

\section{Proof of Proposition \ref{lemma:FM}}
\label{sec:pf-lemma-FM}

In this section, we prove Proposition \ref{lemma:FM}, a key ingredient
in the proof of Theorem \ref{cor:mean} and, in turn, of Theorem
\ref{thm:LA-p_leq_m}. In essence, the proposition asserts that the
sequence $f_A$ from \eqref{eq:f_A} uniquely determines $|A|$ and, up
to internally matched pairs, $A$.  The idea of the proof is that,
treating $f_A$ as an almost periodic sequence in the sense of
Besicovitch, we consider its Fourier expansion. Then $|A|/2$ is the
mean value of $f_A$, and we show that $A$, 
modulo internally matched pairs, is determined by the expansion.  The
latter fact is established by an induction with respect to
divisibility in $S^1$ applied to an explicit expression via $A$ of the
Fourier coefficients.

The proof is structured as follows. In the first part, we recall the
relevant definitions and results from the theory of almost periodic
functions, for they are technically not quite standard, but we refer
to, e.g., \cite{Be, Co} for the proofs. This part is included only for
the sake of completeness and the reader's convenience, and can be
skipped. In the second part, we spell out the Fourier expansion of
$f_A$ in terms of $A$. Finally, in the third part, we show that this
expansion determines $A$ once the internally matched pairs are
excluded.

\begin{proof}[Proof of Proposition \ref{lemma:FM}]
  Without loss of generality, we may assume that $0<\alpha_i<1$ and
  likewise $0<\beta_j<1$ for all $\alpha_i$ and $\beta_j$. As a
  consequence, $\alpha_i$ and $\beta_j$ form an externally matched
  pair if and only if $\alpha_i=\beta_j$, and $\alpha_i$ and
  $\alpha_j$ form an internally matched pair if and only if
  $\alpha_j=1-\alpha_i$.  The fact that $f_A-f_B=\const$ once the
  matching condition is satisfied is an immediate consequence of the
  fact that $\{k\alpha\}+\{k(1-\alpha)\}=1$ for every irrational
  $0<\alpha<1$.  Hence, we need to prove that the condition is
  necessary and calculate the constant.

  \smallskip\noindent \emph{Part 1: Background facts about almost
    periodic sequences.}  For $\lambda\in [0, 1)$, denote by
  $\fe_\lambda$ the exponential sequence
  $\fe_\lambda(k)=\exp\big(2\pi\sqrt{-1}\lambda k\big)$, $k\in\N$. By
  definition, a trigonometric polynomial $P(k)$ is a finite linear
  combination of the sequences $\fe_\lambda$ with complex
  coefficients. Let $\PP$ be the space of all trigonometric
  polynomials on $\N$. The mean value of $P$ is the coefficient of
  $\fe_0$. It can also be evaluated as
  \begin{equation}
    \label{eq:Mean}
    M(P):=\lim_{N\to \infty}\frac{1}{N}\sum_{k=1}^{N} P(k).
  \end{equation}
  The limit indeed exists and is equal to the $\fe_0$-coefficient in
  $P$. Furthermore, consider the inner product on $\PP$ defined as
  $$
  \left<P_1,P_2\right>:=M(P_1\bar{P}_2).
  $$
  One can check that $\|P\|:=\left<P,P\right>^{1/2}$ is a genuine norm
  on $\PP$. More explicitly, for
  $P=\sum_\lambda c_\lambda \fe_\lambda$ (a finite sum), we have
  $$
  \|P\|^2=\sum_\lambda |c_\lambda|^2.
  $$
  The space $\mathfrak{B}_2$ of almost periodic functions in the sense
  of Besicovitch is the completion of $\PP$ with respect to this
  norm. This is a non-separable Hilbert space, and $\fe_\lambda$,
  $\lambda\in [0, 1)$, is an orthonormal basis. The mean $M$ naturally
  extends to a bounded linear functional
  $M\colon \mathfrak{B}_2\to \C$, which is simply the inner product
  with $\fe_0$, i.e., $M(f)=\left<f,\fe_0\right>$.

  To be able to say when a genuine sequence $f$ represents an element
  of $\mathfrak{B}_2$, we need a different definition of this
  space. Let us replace the limit by the upper limit in
  \eqref{eq:Mean} and then in the definition of the norm. Thus for $f$
  we set
   $$
   \|f\|_{\CM}^2= \limsup_{N\to \infty}\ \frac{1}{N} \sum_{k=1}^{N}
   |f(k)|^2.
    $$
    Let $\CM$ be the linear space of $f$ such that
    $\|f\|_{\CM}<\infty$. Then $\|\cdot\|_\CM$ is a semi-norm on
    $\CM$. (This is not a true norm. For instance, $\|f\|_\CM=0$ when
    $f\in \ell^1(\N)$.)  Taking the quotient $\CM/\LL$, where
    $\LL\subset \CM$ is the linear subspace formed by $f$ with
    $\|f\|_\CM=0$, we obtain a Banach space, called the Marcinkiewicz
    space; see \cite{Be, Co}. (It is not obvious from the definition
    that this space is complete with respect to the induced norm.)
    Finally, $\mathfrak{B}_2$ is the image in $\CM/\LL$ of the closure
    of $\PP$ in $\CM$ with respect to $\|\cdot\|_\CM$. Thus, we say
    that $f$ represents an element of $\mathfrak{B}_2$ -- or, for the
    sake of brevity, that $f$ belongs to $\mathfrak{B}_2$ -- whenever
    there is a sequence of trigonometric polynomials $P_n$ converging
    to $f$ in $\|\cdot\|_\CM$, i.e., $\|f-P_n\|_\CM\to 0$. It is worth
    keeping in mind that different sequences $f$ can give rise to the
    same element of $\mathfrak{B}_2$.

    The Fourier expansion of $f\in\mathfrak{B}_2$ is defined as
    follows. For $\lambda\in [0,1)$, set
  $$
  c_\lambda=\left<f,\fe_\lambda\right>.
  $$
  One can show that $c_\lambda=0$ for $\lambda$ outside a finite or
  countable set $\Lambda\subset [0,1)$. The elements of $\Lambda$ are
  called the \emph{Fourier exponents} of $f$, while $c_\lambda$ are
  the \emph{Fourier coefficients}. We write
  $$
  f\sim \sum_{\lambda\in\Lambda} c_\lambda \fe_\lambda
  $$
  and call the right-hand side the \emph{Fourier expansion} of $f$. It
  converges to $f$ in $\mathfrak{B}_2$ and $c_0=M(f)$. For instance,
  when $f=P$ is a trigonometric polynomial, its Fourier expansion is
  $P$ itself. Note, however, that the Fourier series of $f$ need not
  converge pointwise on $\N$ even when $f$ is a genuine sequence.

  \smallskip\noindent \emph{Part 2: The Fourier expansion of $f_A$.}
  Assuming throughout that $\alpha\in (0,1)$ is irrational, set
  $f_\alpha(k)=\{k\alpha\}$, $k\in\N$. It is not hard to see that
  $f_\alpha\in\mathfrak{B}_2$ and $M(f_\alpha)=1/2$. (However,
  $f_{\alpha}$ is not a uniform limit of trigonometric polynomials on
  $\N$. In other words, $f_\alpha$ is not an almost periodic function
  in the sense of Bohr.)  Hence, $f_A\in \mathfrak{B}_2$ and
  $M(f_A)=|A|/2$. As a consequence, $f_A-f_B=(|A|-|B|)/2$ when
  $f_A-f_B$ is constant. Therefore, the condition that
  $f_A-f_B=\const$ is equivalent to
  \begin{equation}
    \label{eq:f_A-f_B}
  f_A-|A|/2=f_B-|B|/2.
\end{equation}

Next, let us modify $A$ by removing all internally matched pairs.
Then \eqref{eq:f_A} turns into
  \begin{equation}
    \label{eq:f_A-modified}
    f_A(k)=\eta+\sum_i \{k\alpha_i\},
  \end{equation}
  where $\eta$ is the number of internally matched pairs in the
  original multiset, and we have now assumed that $A$ contains no such
  pairs. This is a consequence of the fact that
  $\{k\alpha\}+\{k(1-\alpha)\}=1$ for all $k$, and hence every such
  pair contributes one to $f_A$. As before, we also require that
  $0<\alpha_i<1$ for all $\alpha_i\in A$. Note that
  \eqref{eq:f_A-modified} has the same form as the left-hand side of
  \eqref{eq:f_A-f_B}.

  With this in mind, \emph{it suffices to show that $f_A$ given by
    \eqref{eq:f_A-modified} with $\eta\in\R$ completely determines $A$
    whenever $A$ has no internally matched pairs, which we will assume
    from now on}. (As a side remark, note that $f_A$ also determines
  $\eta$. Indeed, when $f_A$ determines $A$, it determines $|A|$ and
  $M(f_A)=\eta+|A|/2$.)  Applying the same argument to $f_B$, we see
  that $A=B$ up to internally matched pairs when $f_A-f_B=\const$.  To
  complete the proof, we will show that $A$ is determined by the
  Fourier expansion of $f_A$ and hence by $f_A$.

  It is not hard to show that the Fourier expansion of $f_\alpha$ is
$$
f_\alpha\sim \frac{1}{2} - \sum_{l \neq 0} \frac{1}{2\pi \sqrt{-1} \,
  l} \ \fe_{l \alpha}.
$$
Hence
$$
f_A\sim \eta +\frac{|A|}{2} -\sum_i\sum_{l \neq 0}
\frac{1}{2\pi\sqrt{-1} \, l} \ \fe_{l \alpha_i}.
$$
In other words, let $\Gamma$ be the multiset of all multiples
$\lambda=l\alpha_i$, for $\alpha_i\in A$ and non-zero integers
$l$. Then the set $\Lambda$ of Fourier exponents of $f_A$ (with zero
excluded) comprises $\lambda\in\Gamma$ such that $c_\lambda\neq 0$.
Here the Fourier coefficient $c_\lambda$ of $\lambda\in\Gamma$ is
$$
c_\lambda=\frac{-1}{2\pi\sqrt{-1}}C_\lambda, \textrm{ where }
C_\lambda=\sum_{l \alpha_i\equiv \lambda}\frac{1}{l}
$$
and the congruence is taken modulo one.  Setting $C_\lambda=0$ when
$\lambda\not\in \Gamma$, we can view
$C\colon \lambda\mapsto C_\lambda$ as a function on $[0,1)$. Now, it
suffices to show that this function uniquely determines the frequency
multiset $A$.

\smallskip\noindent\emph{Part 3: From the Fourier expansion of $f_A$
  to $A$.}  In the rest of the proof, we identify $[0,1)$ with
$S^1=\mathbb{R}/\mathbb{Z}$ and view $A$ as a finite multiset in
$S^1$. For $\alpha\in S^1$, let $d(\alpha)\geq 0$ be the multiplicity
of $\alpha$ in $A$ and set $D(\alpha):=d(\alpha)-d(-\alpha)$. Since
$A$ has no internally matched pairs, i.e., $A\cap(-A)=\emptyset$, at
most one of $d(\alpha)$ and $d(-\alpha)$ is non-zero, and hence
$d(\alpha)=\max\{D(\alpha),0\}$. Thus it suffices to show that $D$ is
determined by the function $C\colon \lambda\mapsto C_\lambda$.

Splitting the sum defining $C_\lambda$ into the terms with $l>0$ and
$l<0$, and writing $l\alpha_i=(-l)(-\alpha_i)$ in the latter, we
obtain
\begin{equation}
  \label{eq:C-D}
  C_\lambda=\sum_{l\geq 1}\frac{1}{l}\sum_{\alpha\colon l\alpha
    =\lambda} D(\alpha).
\end{equation}
The sum is finite for every $\lambda$. Indeed, $D$ is finitely
supported and, since every point $\alpha$ of its support is
irrational, for every $\lambda$ there is at most one $l$ such that
$l\alpha=\lambda$. The right-hand side is linear in $D$. Thus, we can
view \eqref{eq:C-D} as a map from the set of finitely supported
functions $D$ to functions $C$. With this in mind, it is enough to
show that if $D\colon S^1\to\mathbb{Z}$ has finite support $S$
consisting of irrational points and the right-hand side of
\eqref{eq:C-D} vanishes for all $\lambda$, then $D\equiv 0$. In other
words, we need to show that $C\not\equiv 0$ when $D\not\equiv 0$.

Thus, assume that $D\neq 0$, i.e., $S\neq\emptyset$. We call $s\in S$
\emph{primitive} if $s\neq l s'$ for all integers $l\geq 2$ and all
$s'\in S$. A primitive element exists. Indeed, otherwise there would
be a sequence $s_0,s_1,\dotsc$ in $S$ and integers $l_k\geq 2$ such
that $s_{k-1}=l_k s_k$. Since $S$ is finite, $s_j=s_k$ for some $j<k$,
and then $s_j=n s_j$ with $n=l_{j+1}\dotsb l_k\geq 2$. Hence
$(n - 1)s_j=0$ in $S^1$, i.e., $s_j$ is rational, which is
impossible. (To put it differently, divisibility gives rise to a
partial order on $S^1\setminus \Q/\Z$: $s\succ s'$ if and only if
$s=l s'$ with $l\geq 2$. Then the primitive elements are simply the
minimal elements of $S$.) Let $s$ be primitive and set $\lambda=s$ in
\eqref{eq:C-D}. The term with $l=1$ is $D(s)$. All terms with
$l\geq 2$ vanish: $D(\alpha)\neq 0$ only for $\alpha\in S$, and then
$l\alpha\neq \lambda=s$ since $s$ is primitive. Thus $C_s=D(s)\neq 0$,
as desired.

We conclude that $C$ determines $D$ and hence the multiset $A$. This
completes the proof of the proposition.
\end{proof}

\end{document}